\documentclass[11pt]{article}
\usepackage[utf8]{inputenc}
\usepackage{color}
\usepackage{xcolor}
\usepackage{ulem}

\usepackage{amsmath,amsfonts,euscript,amssymb,amsthm,a4,times}

\allowdisplaybreaks
\def\XXint#1#2#3{{\setbox0=\hbox{$#1{#2#3}{\int}$}
\vcenter{\hbox{$#2#3$}}\kern-.5\wd0}}

\newtheorem{theorem}{Theorem}[section]

\newtheorem{lemma}[theorem]{Lemma}

\newtheorem{proposition}[theorem]{Proposition}

\theoremstyle{definition}
\newtheorem{definition}[theorem]{Definition}

\def\tcr{\textcolor{red}} %correzioni
\def\tcb{\textcolor{blue}} %parti nuove
\makeatletter \catcode`@=11
\newbox\tr@tto
\setbox\tr@tto=\hbox{{\count0=0\dimen0=-,9pt\dimen1=1,1pt\loop\ifnum
    \count0<11 \advance \count0 by1 \vrule width.51pt height\dimen1
    depth\dimen0\kern-0.17pt\advance\dimen0 by-0.05pt\advance\dimen1
    by0.1pt\repeat \loop\ifnum\count0<21\advance \count0 by1 \vrule
    width.6pt height\dimen1 depth\dimen0\kern-0.2pt \advance\dimen0
    by-0.1pt\advance\dimen1 by 0.05pt\repeat}}
\def\medint{\displaystyle\copy\tr@tto\kern-10.4pt\int}
\numberwithin{equation}{section}

\def\real{\mathbb R}

\def\proofof#1{\begin{proof}[Proof of #1]}

\def\dx{{\mathrm d}x}
\def\dy{{\mathrm d}y}

\def\loc{\rm loc}

\def\real{\mathbb R}

\def\proofof#1{\begin{proof}[Proof of #1]}

\def\dx{{\mathrm d}x}
\def\dy{{\mathrm d}y}

\def\loc{\rm loc}

\catcode`@=11
\newbox\tr@tto
\setbox\tr@tto=\hbox{{\count0=0\dimen0=-,9pt\dimen1=1,1pt\loop\ifnum\count0<11 \advance \count0 by1 \vrule width.51pt height\dimen1
     depth\dimen0\kern-0.17pt\advance\dimen0 by-0.05pt\advance\dimen1
     by0.1pt\repeat \loop\ifnum\count0<21\advance \count0 by1 \vrule
     width.6pt height\dimen1 depth\dimen0\kern-0.2pt \advance\dimen0
     by-0.1pt\advance\dimen1 by 0.05pt\repeat}}
\def\medint{\displaystyle\copy\tr@tto\kern-10.4pt\int}
\catcode`@=12

\numberwithin{equation}{section}

\date{\today}

\begin{document}

%\maketitle

\title{On the interplay between the anisotropy and the regularity of the coefficients in non autonomous problems}
\author{Michela Eleuteri $\&$ Antonia Passarelli di Napoli}

\maketitle

\begin{abstract}
We establish the local Lipschitz regularity of the local minimizers of non autonomous integral {functionals} of the form
\begin{equation*}
\label{functional}
{\mathcal{F}(z) =} \int_\Omega F(x, Dz)\,dx,
\end{equation*}
where $\Omega$ is a bounded open set of $\mathbb{R}^n$, $n \ge 2$ and the energy density $F(x,\xi)$ satisfies general growth conditions with respect to the gradient variable.
The main novelty of this work lies in the identification of a {suitable}  condition on the interplay between the regularity of the partial map $x\to {F(x,\xi)}$ and the ratio of the non uniform ellipticity of the problem that is sufficient to obtain the Lipschitz regularity of the local minimizers. Our result extends and unifies those previously known from the literature, including the standard {  $p-$growth} case.

\end{abstract}



\bigskip

% \section{To do list}

% - approssimazione seguendo il modello CGM (vedere se funziona)

% - cambiare la stima a priori nel caso vettoriale (?)

% - (Michela:) approssimazione e (Antonella) passaggio al limite

\section{Introduction}
The aim of this paper is to investigate the local Lipschitz regularity of the  local minimizers of variational integrals of the form \begin{equation}\label{functional} \mathcal{F}(z)=\int_{\Omega}F(x,Dz)\,dx, \end{equation} where $\Omega\subset\mathbb{R}^n$, $n\ge2$, is a bounded open set. The energy density $F(x,\xi)$ is assumed to satisfy general growth conditions with respect to the gradient variable. {In the} study of variational problems with nonstandard growth {a key role has been played by} the pioneering works of Marcellini \cite{M89,M91} and, since then, {the regularity theory of functionals with general growth} developed into one of the central topics of modern regularity theory in the Calculus of Variations. The literature on this subject is by now extremely extensive; therefore, rather than attempting a comprehensive account, we refer the interested reader to the survey papers \cite{CM15,M06,Mingione-Radulescu} and the references therein. 
\\ 
When dealing with non-autonomous variational integrals under general growth assumptions. The possible occurrence of the Lavrentiev phenomenon represents one of the main obstacles to the regularity theory. Indeed, the lack of density of smooth competitors may prevent the transfer of the regularity properties from the approximating problems to the minimizers. We refer to the seminal paper \cite{BM} and to the more recent contributions \cite{BDS,FDFL} for a detailed discussion of this phenomenon and its consequences. For this reason, an important line of research has been devoted to identifying sufficient assumptions on the dependence of the energy density upon the spatial variable in order to exclude the Lavrentiev phenomenon and recover the regularity of the minimizers. 
\\
Eleuteri, Marcellini, and Mascolo \cite{EMM16} obtained a fundamental result in this direction. They considered integrands satisfying  $(p,q)$ growth conditions and that, as a function of  the spatial variable, belong to the Sobolev space $W^{1,r}$, with $r>n$. They proved that local minimizers are locally Lipschitz continuous, provided the ellipticity ratio satisfies the sharp closeness condition \[ \frac{q}{p}<1+\frac{1}{n}-\frac{1}{r}. \] This result clearly reveals that the admissible degree of non-uniform ellipticity is directly related to the Sobolev regularity of the coefficients. Later on, these ideas have been considerably extended, leading to a wide regularity theory for vector-valued non-autonomous functionals with non-uniform ellipticity.
However we need to mention that the Lipschitz regularity of minimizers for functionals with general growth in both the autonomous and the non-autonomous case under various regularity assumptions on the coefficients has been a widely studied topic. The list of all contributions in this respect is certainly too wide, therefore we just refer to some recent papers which we believe are closed connected wih the spirit of our contribution, namely 
\cite{BM20, BS20, BB20, BBLV18,   ColMin15, CMMP23,  DFM20, DFM21, DFM23, DFM, DM20, EMM3,EMMP,EPT, ElePerTre, ELM02, GMP24, HO, Mar20, M2021, M2022, MNP26, MP06}. 
Despite the variety of models considered in the above mentioned papers, they all share the same underlying philosophy: one first prescribes the Sobolev regularity of the coefficients (within the Lebesgue scale) and then determines a suitable upper bound on the ellipticity ratio that guarantees the absence of the Lavrentiev phenomenon and, consequently, the regularity of minimizers. From this perspective, the regularity theory is governed by the delicate interplay between two competing features of the energy density: its regularity with respect to the spatial variable and its growth with respect to the gradient variable. 
\\
The main {aim} of this paper is to investigate this interplay further, {by generalizing and complementing the analysis already started in \cite{EMMPdN26} and \cite{brunopini}.} More precisely, our goal is to identify an optimal condition linking the Sobolev regularity of the coefficients to the degree of non-uniform ellipticity of the functional, thereby providing a sharp criterion ensuring the local Lipschitz continuity of minimizers.\\
To  describe our results  precisely, we  now introduce the assumptions.
\\
We shall assume that the energy density $F(x,\xi)$ in \eqref{functional} is {a Carath\'eodory function} such that $\xi \mapsto F(x, \xi)$ is $\mathcal{C}^2(\real^n)$  and   there exists $\tilde F:\Omega\times [0,+\infty)\mapsto {{[0,+\infty)}}$ such that ${\tilde{F}(x, \xi)=F(x,|\xi|)}$. We assume moreover that $x \mapsto {F_{\xi }(x, \xi)}$ is weakly differentiable, with $F_{\xi x}(x, \xi)$ a Carath\'eodory function, and that there exist an exponent $p>1,$  {constants $\lambda,\Lambda>0$, a nonnegative measurable function $h(x)$ and} a non decreasing function $g:[0,+\infty]\to {[0,+\infty]}$ with $g(1)=1$ such that 
\[\lambda \, (1 + |\xi|^2)^{\frac{p-2}{2}} |\mu|^2 \le {{\sum_{i,j  = 1}^n}} {F_{\xi_i \xi_j}}(x, \xi) {\mu_i \mu_j}, \eqno{\rm (F1)}\]
\[
|{F_{\xi_i \xi_j}}(x, \xi)| \le \, \Lambda \, (1 + |\xi|^2)^{\frac{p-2}{2}} g\big((1+|\xi|^2)^{\frac12}\big), \eqno{\rm (F2)} \]
\[|F_{\xi x}(x, \xi)| \le \, h(x) \, (1 + |\xi|^2)^{\frac{p-1}{2}}g\big((1+|\xi|^2)^{\frac12}\big),  \eqno{\rm (F3)}
\]
for every $\xi, \mu \in \mathbb{R}^n$, and a.e. $x$ in $\Omega$. 
Assumptions (F2) e (F3) entail the following non standard growth condition for the energy density $F(x,\xi)$ as a function of the gradient variable 
\[
c(\lambda,p) (1 + |\xi|^2)^{p/2} \le \, F(x, \xi)  \le \, C(\Lambda) (1 + |\xi|^2)^{p/2} g\big((1+|\xi|^2)^{\frac12}\big)
\]
and also implies that the ellipticity ratio of the matrix $a_{ij}={F_{\xi_i \xi_j}}(x, Du)$, that is the ratio between the largest and
the lowest eigenvalue,  is  bounded  by $g(|Du|)$ that may blow up as $|Du|\to\infty$ { (we refer for instance to \cite{DF26} for a unified and comprehensive survey concerning the recent advances of the nonuniformly elliptic setting, with a particular reference to the basic model examples capturing some of the main aspects of the uniform/nonuniform ellipticity notions).}
\iffalse
Hence, we need  to
impose  a bound on the rate of the possible blow-up  of the non-uniform ellipticity of the Euler Lagrange equation of the functional $\mathcal{F}(z)$ that is  for the function $g(t)$ appearing in assumptions {\rm (F2)}, {\rm (F3)}. Actually, we shall assume that there exist a positive strictly increasing function $\mathcal{L}_1(t),$  diverging at $\infty$ and $t_0\ge 0$ such that
\[
g^{{\frac{\tcr{2\cdot}2^*}{(2^* - 2) p}}}(t) \le \, \frac{t}{\mathcal{L}_1(t)}, \qquad \forall t \ge \, t_0, \eqno{\rm (F4)}
\]
where
\[
2^* = \left \{
\begin{array}{lll}
\!\!\!\!\!\! & \displaystyle \frac{2n}{n-2} \qquad & n \ge 3\\[3mm]
\!\!\!\!\!\! & \displaystyle \textnormal{any number greater than 2} \qquad & n = 2.
\end{array}
\right.
\]
\fi
%\colorbox{yellow}{non lo direi qui}
%\textcolor{green}{Hence, by Lemma \tcr{\ref{GGPnew}}, we have that
%\begin{equation}
%\label{conseguenzaLemma2.1g}
%\frac{t}{\mathcal{L}(t)}\lesssim [(\mathcal{L}^{-1})^{*}]^{-1}(t)
%\end{equation}
%which is clearly concave as it is the inverse of a polar function which is always convex.}
%
\\
For the function $h(x)$ appearing in assumption (F3), we shall assume that there exists a convex continuous strictly increasing function $\Phi:[0,+\infty]\to [0,+\infty]$  such that 
%\begin{equation}\label{assumponh}
%\[
% \int_{B_R} h^n \Phi(h)^{n} \, dx < + \infty \qquad \forall B_R \Subset \Omega.  \eqno{\rm (F5)}
%\]
\[
 \int_{B_R} h^n \Phi(h)^{\theta} \, dx < + \infty \qquad \forall B_R \Subset \Omega,  \eqno{{\rm {(F4)}}}
\]
where $\theta$ is a positive exponent such that
\begin{equation}
\label{deftheta}
\left \{
\begin{array}{lll}
\!\!\!\!\!\! & \theta = n \qquad & n \ge  3\\[3mm]
\!\!\!\!\!\! & \theta > 2 + \mu\qquad & n = 2,
\end{array}
\right.
\end{equation}
{for some $0<\mu<\frac14.$}

For the inverse of the function $\Phi(t)$ we moreover assume that
\begin{equation}
\begin{cases} &\log \Phi^{-1}(e^t) \,\,\, \textnormal{is convex} 
\\
& \log \Big(e+\Phi^{-1}(t)\Big)\le c_\Phi t^\mu, \,\,\, \textnormal{for} \,\,\, t>t_0,
\end{cases}\tag{\rm {F5}}
\end{equation}
{where $\mu$ is the same of  assumption \eqref{deftheta}.} 

From {${\rm(F5)}_2$} we deduce that 
\begin{equation}
\label{conseguenzasuPhi}
\Phi(s)\ge \frac{1}{c_\Phi}\log^{\frac{1}{\mu}}(e+s) {\ge \frac{1}{c_{\Phi}}.}
\end{equation}
Observe that assumption {(F4)} implies that, for every fixed $\xi$, the map \( x\longmapsto F(x,\xi) \) belongs to the Orlicz--Sobolev space generated by a Young function with asymptotic growth \[ t^n\Phi^\theta(t) \qquad\text{as }t\to\infty. \]

This class is strictly contained in the critical Sobolev space $W^{1,n}$, consistently with the fact that, even in the case of standard growth (corresponding to $g\equiv\mathrm{const}$ in (F2)--(F3)), Lipschitz regularity is known to hold under the assumption that the coefficients belong to $W^{1,n}\log^\alpha L$, provided that $\alpha$ is sufficiently large depending on the dimension. We refer to \cite{CGHPdN20} and, in the variable exponent setting, to \cite{EPdN23}. 
\\
Here, the interaction between the regularity of the coefficients and the non-uniform ellipticity of the functional is quantified by the auxiliary function \begin{equation}\label{defpsi} \Psi(t):=g^2(t)\Phi^{-1}\!\left(g^2(t)\right). \end{equation} The appearance of the function $\Psi$ is not accidental. Indeed, even in the autonomous case, \[ F(x,\xi)\equiv F(\xi), \] where no regularity of the coefficients is involved, suitable restrictions on the ellipticity ratio are unavoidable. To illustrate this point, recall that for functionals with $(p,q)$-growth the ellipticity ratio behaves asymptotically as \[ |\xi|^{\,q-p}, \] and local Lipschitz continuity is known to hold only under closeness conditions between the exponents. In the general case one requires \[ \frac{q}{p}\le c(n), \qquad c(n)\longrightarrow1 \quad\text{as }n\to\infty, \] whereas dimension-independent estimates are currently available only for {\it a priori} bounded minimizers (see \cite{CKP,CKPPisa,CGGP,EPdN24}), and still rely on the assumption \[ q<p+1. \] In particular, the exponent of the ellipticity ratio $|\xi|^{q-p}$ must remain strictly below one, implying that the ellipticity ratio is eventually concave. This observation suggests that eventual concavity should be regarded as the minimal structural requirement on the ellipticity ratio, independently of any regularity assumption on the coefficients. In our framework the ellipticity ratio is described by \[ g\!\left(\sqrt{1+|\xi|^2}\right), \] so that the previous discussion naturally leads to requiring its asymptotic concavity.

However, in case of non-autonomous functionals, the interaction between the regularity of the coefficients and the non-uniform ellipticity of the functional is described by the function $\Psi$ defined at \eqref{defpsi}; hence we need a regularity of the coefficients sufficient to compensate for the loss produced by $\Phi^{-1}$, preserving the asymptotic concavity of ellipticity ratio. From this perspective, assumption {(F6) below}  can be interpreted as the precise balance between these two competing effects: on the one hand, the blow up of the ellipticity ratio produced by the inverse function $\Phi^{-1}$, and on the other hand, the gain provided by the additional regularity of the coefficients. More precisely, {(F6)} requires that the effective ellipticity measured by $\Psi$ remains strictly below the critical linear growth.
\color{black}
This is expressed by assuming the existence of an invertible increasing function $\mathcal{L}:{[0,+\infty)\to [0,+\infty)}$  such that
\begin{equation}
\label{HpsuLameno1}
\left \{
\begin{array}{lll}
\displaystyle \textnormal{it inverse satisfies {$\mathcal{L}^{-1}(0) = 0$}}\\[2mm]\textnormal{{the map}
$\displaystyle t \mapsto \frac{\mathcal{L}^{-1}(t)}{t}$
is positive and increasing on $[t_0, + \infty)$ with $t_0 > 0$}\\[2mm]
\textnormal{{the map}
$\displaystyle t \mapsto \frac{\mathcal{L}^{-1}(t)}{t}$  verifies 
$\displaystyle
\lim_{t \rightarrow 0} \frac{\mathcal{L}^{-1}(t)}{t} = 0, \quad \lim_{t \rightarrow + \infty} \frac{\mathcal{L}^{-1}(t)}{t} = + \infty;$}
\end{array}
\right.
\end{equation}
\color{black}
{moreover $\mathcal{L}$ is} such that \begin{equation*}\label{assumptionpsi} \Psi^{\frac{2^*}{(2^*-2)p}}(t) = \left( g^2(t)\Phi^{-1}\!\left(g^2(t)\right) \right)^{\frac{2^*}{(2^*-2)p}} \le \frac{t}{\mathcal L(t)}, \qquad t\ge t_0, \eqno{{\rm(F6)}} \end{equation*} 
{where we recall that
\[
2^* = \left \{
\begin{array}{lll}
\!\!\!\!\!\! & \displaystyle \frac{2n}{n-2} \qquad & n \ge 3\\[3mm]
\!\!\!\!\!\! & \displaystyle \textnormal{any number greater than 2} \qquad & n = 2.
\end{array}
\right.
\]
From now on, for the sake of simplicity, without loss of generality, we shall assume $t_0 = 1.$}

Note that by {(F6)} and the monotonicity of $\Phi^{-1}$, since $g(s)\ge 1$ for $s\ge 1$, we deduce that
\begin{eqnarray*}  &&\Big( g^2({t})\Phi^{-1}\!(1) \Big)^{\frac{2^*}{(2^*-2)p}} \le  \Big( g^2({t})\Phi^{-1}\!(g^2(1))\Big)^{\frac{2^*}{(2^*-2)p}}\cr\cr &\le&  \Big( g^2({t})\Phi^{-1}\!(g^2({t})) \Big)^{\frac{2^*}{(2^*-2)p}} \le  \frac{{t}}{\mathcal L({t})}, {\qquad t\ge 1.} \qquad  \end{eqnarray*}
This implies, in particular, that
\begin{equation}
\label{conseguenzasug}
g^{\frac{2\cdot2^*}{(2^*-2)p}}({t})\le C \frac{{t}}{\mathcal L({t})}, {\qquad t\ge 1.}
\end{equation}

Therefore, assumption {(F6)} should not be viewed merely as a technical hypothesis. Rather, it identifies the threshold at which the additional regularity of the coefficients exactly compensates for the non-uniform ellipticity of the functional, revealing the optimal balance between these two competing mechanisms.
%Indeed  in this case by Lemma \tcr{\ref{GGPnew}}, we have that
%\begin{equation}
%\label{conseguenzaLemma2.1Psi}
%\frac{t}{L(t)}\lesssim [(L^{-1})^{*}]^{-1}(t)
%\end{equation}
%\colorbox{yellow}{mettiamo qui anche (1.2)?}
%which is clearly concave as it is the inverse of a polar function which is always convex.
%\tcr{Condition (F8) obviously is equivalent to require that
%$$\Psi(t) \le \, \frac{t^{\textcolor{magenta}{\frac{2^*}{(2^* - 2) p}}}}{{L}^{{\textcolor{magenta}{\frac{2^*}{(2^* - 2) p}}}}(t)} \qquad  t\ge t_0   $$}
%\\
%\colorbox{yellow}{controllare questa stima!}
\\
%\colorbox{yellow}{Ho aggiunto ipotesi su $\mathcal{L}$ in modo da poter applicare il Lemma 3.6}

%\begin{equation*}
% \Psi^{{{\frac{2^*}{(2^* - 2) p}}}}(t)+g^{{{\frac{\tcr{2\cdot}2^*}{(2^* - 2) p}}}}(t) \le \, \frac{t}{\mathcal{L}(t)}\qquad  t\ge t_0. \eqno{\rm (F9)}
%\end{equation*} 

We now recall the definition of local minimizer
\begin{definition}
A mapping $u\in W^{1,1}_{\rm loc}(\Omega)$ is a local minimizer {of the integral functional \eqref{functional}} if
$F(x,Du) \in L^{1}_{\rm loc}(\Omega )$ and
\begin{equation}\label{minineq}
 \int_{\mathrm{supp}\varphi} \! F(x,Du) \, dx \leq \int_{\mathrm{supp}\varphi} \! F(x,Du+D\varphi) \, dx  
\end{equation}
for  any $\varphi\in W^{1,1}(\Omega)$ with $\mathrm{supp}\varphi\Subset\Omega$.
\end{definition} 
Actually, we are able to prove the following
\begin{theorem}
    Let $u\in W^{1,p}_{\mathrm{loc}}(\Omega)$ be a local minimizer of \eqref{functional} under the assumptions {{\rm (F1)--(F6)}}. %If the \tcm{function $\Psi$}  %$$\Psi(t)=g^2(t)\Phi^{-1}(g^2(t))$$ 
 %\tcm{introduced in \eqref{defpsi} satisfies {\rm (F6)}}, 
 {Then}
 $u\in W^{1,\infty}_{\mathrm{loc}}(\Omega)$ and the following estimate
\begin{eqnarray}\label{estfin}
	||Du||_{L^{\infty}(B_{\rho})} &\le&   \left(\mathcal{L}^{-1}\right)^{**}\Bigg(C\left(\int_{B_{ R}}(1+ |Du|)^{p}
\, \dx\right)^{\frac{1}{p}}\Bigg)
\cr\cr
&&+\frac{C}{(  R- \rho )^{\frac{4\cdot2^*}{p(2^*-2)}}}\left(\int_{B_{  R}}(1+ |Du|)^{p}
\, \dx\right)^{\frac{1}{p}}
%&& {+ \Big(\tilde g^{-1}\Big)^*\Bigg(\frac{C}{(R-\rho )^{\frac{4\cdot2^*}{p(2^*-2)}}}\left(\int_{B_{{R}}} (1+ |Du|)^{p}
%\, dx\right)^{\frac{1}{p}}\Bigg),}
\end{eqnarray}
 holds true for every $B_\rho\subset B_R\Subset\Omega,$ with a constant $C=C(2^*,p,\lambda,\Lambda, h,\Phi,\theta,\mu)$. 
Here  $\left(\mathcal{L}^{-1}\right)^{**}(\tau)$ denotes the bipolar conjugate of the function $\mathcal{L}^{-1}(\tau)$, where $\mathcal{L}(t)$ is the function appearing in assumption {{\rm {(F6)}}}.
% {Here} we set $\tilde \Psi(t)=\Psi^{\frac{n}{2p}}(t),$ {being $(\tilde{\Psi}^{-1})^*$ the polar conjugate of the function $\tilde{\Psi}^{-1}$} {and $\tilde g(t)= g^{\frac{n}{2p}}(t),$ {being $(\tilde{g}^{-1})^*$ the polar conjugate of the function $\tilde{g}^{-1}$}}.
\end{theorem}

The proof of the previous theorem proceeds in three steps.
\begin{itemize}
\item First, an \textit{a priori} estimate is proved for \textit{smooth minimizers} of $\mathcal{F}(u)$, with constants depending only on the structural constants in assumptions {(F1)--(F6)};
\vskip.2cm
\item Next, the functional in \eqref{functional} is approximated by a sequence of regular functionals $\mathcal{F}_h(u)$, each admitting a unique \textit{smooth} minimizer $u_h$;
\vskip.2cm
\item Finally, the \textit{a priori} estimate is applied to the sequence $(u_h)$, and it is shown that $u_h$ converges to a minimizer of the original functional, preserving the corresponding regularity properties.
\end{itemize}
The main novelty of the present paper lies in the first step of the proof, where the crucial issue is to identify the optimal balance between the regularity of the coefficients and the growth of the ellipticity ratio. More precisely, the additional regularity of the coefficients beyond the critical space $W^{1,n}$ must compensate for the non-uniform ellipticity by sufficiently slowing down the asymptotic growth of the ellipticity ratio. Once this fundamental {a priori} estimate has been established, the remaining part of the proof follows a rather standard approximation scheme. Indeed, thanks to the Uhlenbeck structure of the integrands, one can construct a family of regularized functionals whose smooth minimizers inherit the desired estimates, and then pass to the limit by means of standard compactness arguments.

%\color{blue}
%It worth pointing out that   most of the representative literature cases can be described in terms  of regular varying functions that provide a unified  framework for a broad class of Orlicz growth functions with power-like behavior, including slowly varying perturbations. 
%Hence, we conclude this introduction giving a sufficient condition for our assumption (F9) to hold in this context.

\color{black}

\section{{Comparison with representative literature cases}}

\label{sectionesempi}

This section is devoted to a comparison with the existing results in the literature. In particular, we show that our theorems cover and extend the known results, including those concerning the case of standard growth conditions, that is the {\it  $p-$growth case} as well as those concerning autonomous problems.
\begin{enumerate}
    \item  {\bf Standard growth conditions.}  The standard $p$-growth case corresponds to assume  $g(t)=constant$ in {(F2)--(F3).}
    %\textcolor{magenta}{In this case Proposition \ref{propPsieg} is trivially satisfied as long as both $g \in RV_0$ and $\Psi \in RV_0.$}
    Assumption {(F6)} is clearly satisfied for every strictly increasing function $\mathcal{L}(t)$, while assumption {(F5)$_2$} yields that $\Phi(t)\ge c \log^{\alpha}t,$ with $\alpha>4$ and for $t\ge t_0$. 
    Therefore, a sufficient condition to the Lipschitz regularity of the minimizers is that the partial map $x\to F_\xi(x,\xi)$ belongs to the Orlicz-Sobolev space associated with the Young function 
    %\[t \mapsto t^n \Phi^{\theta}(t)\]
   % where $\Phi \in RV_0$ is a positive, slowly varying function satisfying 
    %$\Phi(t)\ge c \log^{4}t$ , for $t\ge t_0$. The classical logarithmic assumption is recovered by choosing 
    \[
    \Phi(t) = \log^{\alpha}(e + t)
    \]
    with $\alpha > 4.$ This recovers the result in \cite{CGHPdN20} and  can be compared with the corresponding Lipschitz regularity results in \cite{EPdN23} where the variable exponent case is treated under the same assumptions. This is not surprising, since although variable growth falls within the class of non-standard growth conditions, the methods used to establish regularity in this setting are quite similar to those employed in the standard growth case.

    \item {\bf $(p,q)-$growth conditions.}
{Let us consider, for simplicity,  $ n \ge 3$ .} This case can be retrieved by choosing in assumptions {(F2)--(F3),}  $g(t) = t^{q - p}$. In order to have {\eqref{conseguenzasug}} satisfied we need to impose  
$$g^{\frac{n}{p}}(t)=t^{\frac{n(q-p)}{p}}\le \frac{t}{{\mathcal{L}(t)}},$$
for a strictly increasing function {$\mathcal{L}(t)$ satisfying \eqref{HpsuLameno1}}. This implies
%In this case  assumption at (F6) on $\Psi$ reads as $$ \tcr{\Psi}(t)=\left[t^{2(q-p)}\Phi^{-1}(t^{2(q-p)})\right]^{\frac{n}{2p}} {\le \, \frac{t}{\mathcal{L}(t)}} \qquad {\forall t \ge t_0}$$
%{where $\mathcal{L}(t)$ is a suitable increasing positive function, and it can be satisfied by choosing}
%$$\Phi(t)=t^{\frac{n(q-p)}{p-(q-p)n}} {{L}(t)},$$ for some increasing function ${{L}(t)}$.
%Indeed, one can easily check that
%$$\tcr{\Psi}(t) \, {\cong} \left[t^{2(q-p)}\cdot \frac{t^{2(q-p)\frac{p-(q-p)n}{n(q-p)}}}{ L(t)}\right]^{\frac{n}{2p}} {\cong} \frac{t}{\mathcal{L}(t)},$$
%with \tcr{$\mathcal{L}(t)=\left(L(t)\right)^{\frac{n}{2p}}$} . 
 %Note that $\Phi(t)$ is increasing iff 
 
 $$p-(q-p)n>0\,\,\Longleftrightarrow\,\, q<p+\frac{p}{n}.$$ 
 For every $q$ in this range, to satisfy {(F6)} it suffices to choose $$\Phi(t)\approx t^{\frac{n(q-p)}{p-n(q-p)}}L(t),$$ 
with a strictly increasing function $L(t)$ whose inverse satisfies \eqref{HpsuLameno1}. Indeed, with such a choice
 $$\Phi^{-1}(s)\approx \frac{s^{\frac{p-n(q-p)}{n(q-p)}}}{L^{\frac{p-n(q-p)}{n(q-p)}}(\Phi^{-1}(s))}$$ and {(F6)} reads as  $$ t^{\frac{n(q-p)}{p}}\left(\Phi^{-1}(t^{2(q-p)})\right)^{\frac{n}{2p}}=\frac{t}{L^{\frac{p-n(q-p)}{2p(q-p)}}(\Phi^{-1}(t^{2(q-p)}))}=: \frac{t}{\mathcal{L}(t)},  $$
 and it is clearly satisfied if {$L(t)$ is a strictly increasing function satisfying \eqref{HpsuLameno1}}.
 Therefore  the Lipschitz continuity of the minimizers holds true provided the function $h(x)$ in assumption (F4) satisfies
$$\int_{B_R}h^{n\left(1+\frac{n(q-p)}{p-n(q-p)}\right)}L^n(h)\,dx=\int_{B_R}h^{\frac{np}{p-n(q-p)}}L^n(h)\,dx<+\infty.$$ 
Hence, the minimizers of a functional satisfying (F1)--(F3) with $q$ such that
$$p\le q<p+\frac{p}{n}$$
are Lipschitz continuous provided {(F4)} is satisfied with a function $h(x)$ belonging to a Sobolev space {strictly} contained in $$W^{1,\frac{np}{p-n(q-p)}}_{\mathrm{loc}}(\Omega)$$ and a good choice is to assume that the map $x\to F_\xi(x,\xi)$ belongs to  the space
$$W^{1,\frac{np}{p-n(q-p)}\log^\alpha L}_{\mathrm{loc}}(\Omega), \qquad \alpha>0.$$
Note that, if we set \[
{r = \frac{np}{p - n (q-p)}},
\] we have $$r=\frac{np}{p-n(q-p)}\,\, \Longleftrightarrow\,\, \frac{q}{p}=1+\frac1n-\frac1r$$
{hence we retrieve Theorems 1.2 and 1.3 in \cite{EMMPdN26}. On the other hand, if 
\[
r > \frac{np}{p - n(q-p)},
\]
we can choose} $\Phi(t) = t^{\frac{r - n}{n}}$ in {(F4)}. Note that with this choice {(F5)} is trivially satisfied  while {(F4)} is equivalent to assume $h \in L^r$. Indeed 
\[
h^n \Phi(h)^{n} = h^n h^{ \frac{r - n}{n} {n}} = h^r.
\]
Moreover, {in this case}, we  have %(essentially)
\[
\Psi(t) \, {\cong} \, t^{2(q - p) + \frac{2n}{r - n}(q-p)}
\]
%\textcolor{magenta}{that assumtion (F9) is satisfied due to the fact that Proposition \ref{propcasosub} holds with the choices 
%\[\sigma = q - p \qquad \textnormal{and} \qquad a = \frac{r - n}{n}.\]
%In this case indeed \eqref{HpPhieg} is satisfied because}
%\[2 (q - p) \left (1 + \frac{n}{r - n} \right ) < \frac{2p}{n}\]
%i.e.
and assumption {(F6)} is satisfied since
\begin{equation}
\label{cond-EMM}
q < p \left (1 + \frac{1}{n} - \frac{1}{r} \right ).
\end{equation}
%and this is exactly the requirement in the result by Eleuteri-Marcellini-Mascolo \cite{EMM16} {(see also Theorem 1.4 in \cite{EMMPdN26}).}

We would like also to point out that our method allows to extends the same obtained by Marcellini in \cite{M91} with Lipschitz continuous coefficients, also in the case of weaker assumptions on the partial map $x\to F_\xi(x,\xi)$. 
More precisely, the same sharp bound found by Marcellini in \cite{M91} is  sufficient for the Lipschitz regularity of the minimizers also in case of coefficients whose derivative are exponentially integrable. Indeed, for every $p\le q<p+\frac{p}{n}$, we could  choose $\Phi(t) = \exp^{\beta}(t)$ in {(F4)} so that
\[
[\Psi(t)]^{\frac{n}{2p}} = [2^{\frac{1}{\beta}} g^2(t) \log^{\frac{1}{\beta}} g(t)]^{\frac{n}{2p}} \approx t^{\frac{q - p}{p} n} \log^{\frac{n}{2\beta p}}(t)
\]
and so (F6), i.e.
$$\Psi^{\frac{n}{2p}}(t)=t^{\frac{q - p}{p} n} \log^{\frac{n}{2\beta p}}(t)\le \frac{t}{\mathcal{L}(t)}$$
holds true choosing  $$\mathcal{L}(t)=\frac{t^{1-\frac{q - p}{p} n}}{\log^{\frac{n}{2 \beta p}}(t)}$$ which is strictly increasing if %is satisfied if
\[
\frac{q - p}{p} n < 1 \qquad \Leftrightarrow \qquad q < p \left (1 + \frac{1}{n} \right ).
\]
This  turns to be the limit condition in \eqref{cond-EMM} as $r \rightarrow \infty,$ consistently with the fact that  every function {with weak derivatives exponentially integrable locally in $\Omega$}   belongs to $W^{1,r}_{\mathrm{loc}}(\Omega),$  for every $ r > n.$

\item {\bf Case where the ellipticity ratio has slow growth}

When the ellipticity ratio $g(t)$ has slow growth at infinity, namely when we deal with a functional with a mild anisotropy, for example $g(t) = \log^{\beta} t$, with $\beta>0$, then it is sufficient to assume $\Phi(t) =\log^{\alpha} (1+t)$ with {$\alpha > 2 \beta.$}

 Indeed in this case 
%\[g \in RV_0 \Rightarrow g^{\beta} \in RV_0.\] inverse can be calculated exlicitly, by getting
\[
\Phi^{-1}(s) \approx \exp (s^{1/\alpha}),
\]
thus
\[
\Psi(t) \approx \log^{2\beta} t \exp (\log^{\frac{2\beta}{\alpha}}(1+t)).
\]
thus {(F6)} is satisfied, as long as $\alpha > 2 \beta.$ 
%then
%\[0 < \delta := \frac{2 \beta}{\alpha} < 1\]
%and this is a classical example of slowly varying function (see for instance \cite[Section 1.3, Subsection3]{BGT87}.\\
%Thus also $\Psi \in RV_0$ and assumption (F9) is trivially satisfied by Proposition \ref{propRV}.
%}

\item{\bf {Autonomous functionals}}

In case of autonomous functionals, i.e. if $F(x,\xi)\equiv F(\xi)$, assumption (F3) holds true with $h(x)\equiv 0$ and so  {(F4)} is trivially satisfied for every choice of $\Phi(t)$. Therefore, in this case, the Lipschitz regularity of the minimizers holds true under the only assumption (F6) which reduces merely to an assumption on $g.$ To compare this situation with the familiar case of the $(p,q)$ growth conditions, i.e. $g(t)\eqsim t^{q-p}$ at infinity, {condition \eqref{conseguenzasug}} corresponds to assume that
$$ t^{\frac{n(q-p)}{2p}}\le \frac{t}{\mathcal{L}(t)}$$
that is obviously true for every $q$ such that
$$q<p+\frac{p}{2n}.$$
Indeed, in this case, setting $q-p=\vartheta \frac{p}{2n}$, with $0<\vartheta<1$, assumption {(F6)} is satisfied with
$$\mathcal{L}(t)= t^{1-\vartheta}$$
which is strictly increasing. Hence, we recover the bound in \cite{ELM1,ELM02}.
On the other hand, our main result covers in same sense also the limit case $q = \, p + \frac{p}{2n}$, if the function $g(t)$ grows at $\infty$ less than $t^{p + \frac{p}{2n}}$ not necessarily in the scale of the power functions, but for example assuming
\[
g(t) = \frac{t^{p + \frac{p}{2n}}}{L(t)},
\]
for a strictly increasing function $L(t)$.

%\colorbox{yellow}{becchiamo il caso $q \le \, p + \frac{p}{2n}$} con l'uguale se 
%\[
%g(t) = \frac{t^{p + \frac{p}{2n}}}{L(t)}
%\]
%con $L$ crescente

\end{enumerate}
  
%Indeed, from \eqref{legamePsigPhi}, we would have (essentially)
% \[
% \Psi(t^2) \sim \textcolor{red}{constant}
% \]
% therefore inequality \eqref{stimaquasifin} becomes
% \begin{eqnarray}
% 	\textcolor{red}{||V(Du)||_{L^{\infty}(B_{\bar\rho})} \le  \left(\frac{C}{(\bar R-\bar \rho )^{2}} \big(\|V(Du)\|^2_{\infty}\big)\right)^{\frac{n}{p}} \left(\int_{B_{\bar R}}V(Du)^{p}
% \, dx\right)^{\frac{1}{p}},}
% \end{eqnarray}
% and this term in the right hand side can be absorbed only if $2n/p < 1$ which of course can't be realized.

\section{Preliminary {results}}
\medskip
Throughout the paper, the symbol $C$ denotes a generic positive constant, whose value can change from line to line, and even within the same chain of estimates. Whenever necessary, its dependence on the relevant parameters will be explicitly indicated by parentheses or subscripts. For $x\in\mathbb{R}^n$ and $r>0$, we denote by \[ B(x,r)=B_r(x):=\{y\in\mathbb{R}^n:\ |y-x|<r\} \] the open ball centered at $x$ with radius $r$. When no ambiguity arises, we simply write $B_r$ \sout{or $B$}.
\\
We recall the following standard iteration lemma, whose proof can be found, for example, in \cite[Lemma~6.1, p.~191]{Giusti}.
  \begin{lemma}\label{lem:Giaq}
    For $0<R_1<R_2$, consider a bounded function
    $f:[R_1,R_2]\to[0,\infty)$ with
    \begin{equation*}
      f(r_1)\le\vartheta f(r_2)+\frac A{(r_2-r_1)^\alpha}+ \frac B{(r_2-r_1)^\beta}+ C
       \qquad\mbox{for all }R_1\le r_1<r_2\le R_2,
    \end{equation*}
    where $A,B,C$ {denote nonnegative constants, $\alpha > \beta > 0$ 
    and $\vartheta\in[0,1)$.} Then we have
    \begin{equation*}
      f(R_1)\le c({\alpha, \vartheta})
      \bigg(\frac A{(R_2-R_1)^\alpha}+\frac B{(R_2-R_1)^\beta}+C
      \bigg).
    \end{equation*}
  \end{lemma}
  We shall use the elementary inequality contained in the following 
\begin{lemma}\label{elem1} Let $\varphi: (0,+\infty)\to (0,+\infty)$ be  such that the function $\log\,\varphi(e^t)$ is convex. Then
$$\varphi(ab)\le \sqrt{\varphi(a^2)\varphi(b^2)}$$
for every {$a,b > 0$}.
\end{lemma}
\begin{proof}  It suffices to set $$a=e^s\qquad b=e^t$$ and observe that  $$\varphi(ab)=\varphi(e^{s+t})=\varphi(e^{\frac{2s+2t}{2}}).$$ Hence,  the convexity of the function   $\log(\varphi(e^t))$ yields
$$\log\varphi(ab)\le \log\varphi(e^{\frac{2s+2t}{2}})\le \frac12\log\varphi(e^{2s})+\frac12\log \varphi (e^{2t})=\frac12\log\varphi(a^2)+\frac12\log \varphi (b^2)$$
and so
\begin{equation*}  
\log\varphi(ab)\le \frac12\log\big(\varphi(a^2)\varphi (b^2)\big)=\log\big(\sqrt{\varphi(a^2)\varphi (b^2)}\big),
\end{equation*}
i.e. 
the conclusion.\end{proof}
%\tcm{This lemma will be applied twice in the proof of the a priori estimate to the function $\Phi^{-1}$ which satisfies the assumptions in view of (F5)$_2$.} 
%\colorbox{yellow}{dobbiamo proprio dirlo??}
Note that the log-convexity assumption on $\varphi$ in previous lemma is satisfied by most elementary functions, for example
$$\varphi(t)=t^p,\qquad p\in\mathbb{R}$$
$$\varphi(t)=\exp(\log^\alpha(e+t)),\qquad \alpha\ge 1$$
or more generally
$$\varphi(t)=\exp(\psi(\log(e+t))),\qquad \psi\,\,\text{convex}.$$
%\noindent\colorbox{pink}{se non mettiamo l' ipotesi di log-convessit\'a aumentano gli esponenti e non ritroviamo i casi di letteratura.}
%\\
%\\
%\\
%\noindent\colorbox{pink}{ho inserito la sottosezione seguente}
\subsection{An interlude on Orlicz spaces}
Here we introduce a
few auxiliary functions and some basic properties of {Young} functions, that will be used in what follows. 

A function {$\phi:[0,\infty)\to [0,\infty)$ is called a
Young function if and only if it is a left-continuous convex function  that vanishes at $0$ and is not constant. 
It is said to
satisfy the $\Delta_2$-condition if and only if there is a constant
$c>1$ such that $\phi(2t)\le c\,\phi(t)$ for every $t>0$. The smallest
constant $c$ with this property is usually denoted by $\Delta_2(\phi)$.

The conjugate of a {Young} function~$\phi$,  defined for $t \geq 0$ as
\begin{align}
\label{conjugate}
  \phi^*(t) := \sup_{s \geq 0} \{ st - \phi(s) \},
\end{align}
 is  also a Young function and the following Young type inequality
\begin{equation}
  \label{Young-psi}
  st\le \phi(s)+ \phi^*(t)
\end{equation}
holds  for $s,t\ge0$.
\iffalse
For further needs, we recall that for a non-trivial function $\phi:[0,+\infty)\to [0,+\infty)$,
its polar (or Fenchel conjugate)  is defined by
\begin{equation}
\label{Fconiug}
\tcr{\phi}^*(s) := \sup_{t \ge 0} \big \{st - \tcr{\phi}(t) \big \} \qquad {\forall \,}s\ge 0.
\end{equation}
The function $\tcr{\phi}^*: [0,+\infty) \rightarrow [0,+\infty)$ is convex and   the Young-type (or Fenchel) inequality   directly follows from its  definition
\begin{equation}
\label{Fenchel}
st \le \, \tcr{\phi}^*(s) + \tcr{\phi}(t)
\end{equation}
for all $s,t\ge 0$.
\fi
For further needs, we explicitely recall that if  $$\phi (t)\approx t \log^\beta (1+t)\quad \text{near infinity}$$ and $\beta >0$, then 
$$\phi^* (t) \approx \exp (t^{\frac{1}{\beta}})-1\quad \text{near infinity.}$$
The  Orlicz space $L^\phi(\Omega)$, associated with the function $\phi(t)$, is defined as 
\begin{align*}
    \label{orlicz}
    L^\phi (\Omega)= \bigg\{u: \Omega \to \mathbb R^n: \exists \lambda >0\,\,\text{s.t.} \,\, \int_\Omega\phi\bigg(\frac{|u|}{\lambda}\bigg)\, \dx <\infty\bigg\}
\end{align*}
and it becomes a Banach space with the Luxemburg norm
\begin{equation}
\label{Luxemburgnorm}
  \|u\|_{L^\phi(\Omega)}:=\inf\bigg\{\lambda>0\colon \int_\Omega
  \phi\bigg(\frac{|u(x)|}{\lambda}\bigg)\dx\le 1\bigg\}.
\end{equation}
We shall use the common notations
$L\log ^\beta L (\Omega)$ and ${\rm EXP}_\alpha(\Omega)$ to denote the Orlicz spaces associated with the functions $t \log^\beta (1+t)$ and $\exp (t^{\alpha })-1$, respectively.
\\
In the sequel we will use the following equality
\begin{equation}
\label{galquadrato}
\|\phi^2\|_{{\rm EXP}_{\alpha}(\Omega)} = \|\phi\|^2_{{\rm EXP}_{2 \alpha}(\Omega)}
\end{equation}
that easily follows by the definition of the Luxemburg norm at \eqref{Luxemburgnorm}. Indeed, we have 
\[
\|\phi^2\|_{{\rm EXP}_{\alpha}(\Omega)}
=
\inf\left\{\lambda>0:
\int_{\Omega}
\left(
e^{\left(\frac{|\phi|^2}{\lambda}\right)^{\alpha}}
-1
\right)\,\dx
\le 1
\right\}.
\]
Since
\[
\left(\frac{|\phi|^2}{\lambda}\right)^{\alpha}
=
\left(\frac{|\phi|}{\sqrt{\lambda}}\right)^{2\alpha},
\]
setting $\mu=\sqrt{\lambda}$ yields
\[
\begin{aligned}
\|\phi^2\|_{{\rm EXP}_{\alpha}(\Omega)}
&=
\inf\left\{\mu^2>0:
\int_{\Omega}
\left(
e^{\left(\frac{|\phi|}{\mu}\right)^{2\alpha}}
-1
\right)\,\dx
\le 1
\right\} \\
&=
\left(
\inf\left\{\mu>0:
\int_{\Omega}
\left(
e^{\left(\frac{|\phi|}{\mu}\right)^{2\alpha}}
-1
\right)\,\dx
\le 1
\right\}
\right)^2 \\
&=
\|\phi\|_{{\rm EXP}_{2\alpha}(\Omega)}^2.
\end{aligned}
\]
\color{black}
\\
For future reference,  we observe that, if $|\Omega|<\infty$, the following inequalities hold
{\begin{equation*}
  ||f||_{L\log L(\Omega)}\ge 1\,\,\Longrightarrow \,\,||f||_{L\log L(\Omega)}\le  \int_{\Omega}|f|\log\left(e+|f|\right)\,\dx 
 \end{equation*}}
 and for every $\alpha>1$ there exists $c(\alpha)>0$ such that
{
 \begin{equation*}
 ||f||_{L\log L(\Omega)}< 1\,\,\Longrightarrow\,\, ||f||^\alpha_{L\log L(\Omega)}\le C(\alpha)  \int_{\Omega}|f|\log\left(e+|f|\right)\,\dx. 
 \end{equation*}}
 Hence, joining the two cases
 {\begin{equation}\label{alphallogl}||f||_{L\log L(\Omega)}\le \int_{\Omega}|f|\log\left(e+|f|\right)\,\dx +C(\alpha,|\Omega|)\left(\int_{\Omega}|f|\log\left(e+|f|\right)\,\dx\right)^{\frac{1}{\alpha}}\end{equation}}
 Indeed, using again the definition of Luxemburg norm and the fact that $\phi$ satisfies a $\Delta_2$ condition yields {(see \cite{RaoRen91}, Chapter 3, Par. 3.4, Proposition 6)} 
$$\int_{\Omega}\phi\left(\frac{|f|}{||f||_{L^\phi(\Omega)}}\right)\,\dx=1,$$
 for $\phi(t)=t \log (1+t)$  we have
 \begin{equation}
     1=\int_{\Omega}\frac{|f|}{||f||_{L\log L(\Omega)}}\log\left(1+\frac{|f|}{||f||_{L\log L(\Omega)}}\right)\,\dx.
 \end{equation}
 Hence, if $||f||_{L\log L(\Omega)}\ge 1$, we immediately  have
 \begin{equation*}
     1=\int_{\Omega}\frac{|f|}{||f||_{L\log L(\Omega)}}\log\left(1+\frac{|f|}{||f||_{L\log L(\Omega)}}\right)\,\dx\le \int_{\Omega}\frac{|f|}{||f||_{L\log L(\Omega)}}\log\left(1+|f|\right)\,\dx
 \end{equation*}
 and so
 \begin{equation*}
     1\le \frac{1}{||f||_{L\log L(\Omega)}}\int_{\Omega}|f|\log\left(1+|f|\right)\,\dx
 \end{equation*}
 that implies
 \begin{equation*}
     ||f||_{L\log L(\Omega)}\le \int_{\Omega}|f|\log\left(1+|f|\right)\,\dx.
 \end{equation*}
 On the other hand, if $||f||_{L\log L(\Omega)}< 1$
\begin{eqnarray*}
 1&=&\int_{\Omega}\frac{|f|}{||f||_{L\log L(\Omega)}}\log\left(1+\frac{|f|}{||f||_{L\log L(\Omega)}}\right)\,\dx
 \cr\cr
&\le &\int_{\Omega}\frac{|f|}{||f||_{L\log L(\Omega)}}\log\left(\frac{1}{||f||_{L\log L(\Omega)}}+\frac{|f|}{||f||_{L\log L(\Omega)}}\right)\,\dx
 \cr\cr
 &=& \frac{1}{||f||_{L\log L(\Omega)}}\int_{\Omega}|f|\left[\log\left(1+|f|\right)+\log\left(\frac{1}{||f||_{L\log L(\Omega)}}\right)\right]\,\dx\cr\cr
&\le& {\frac{1}{||f||_{L\log L(\Omega)}}\int_{\Omega}|f| \log\left(e+|f|\right) \dx + \log\left(e + \frac{1}{||f||_{L\log L(\Omega)}}\right) \int_{\Omega} |f| \,\dx}\cr\cr
 &\le& \frac{{2}}{||f||_{L\log L(\Omega)}}\log\left(e+\frac{1}{||f||_{L\log L(\Omega)}}\right)\int_{\Omega}|f|\log\left(e+|f|\right)\,\dx\cr\cr
 &\le& {C(\alpha)}\left(\frac{1}{||f||_{L\log L(\Omega)}}+\frac{1}{||f||^\alpha_{L\log L(\Omega)}}\right)\int_{\Omega}|f|\log\left(e+|f|\right)\,\dx\cr\cr
 &\le& 2 {C(\alpha)}\frac{1}{||f||^\alpha_{L\log L(\Omega)}}\int_{\Omega}|f|\log\left(e+|f|\right)\,\dx ,   
\end{eqnarray*}
where, in the second last line, we used that for every  $\alpha>1$ there exists $C(\alpha)>0$ such that $x\log(e+x)\le C(\alpha) (x+x^\alpha)$ for every $x>1$ and in the last line we used that $x<x^\alpha$, for $x>1$ and $\alpha>1$. From the previous estimate \eqref{alphallogl} immediately follows.
\\
\\
The following H\"older type inequality in Orlicz spaces will be useful in what follows.
\begin{proposition}\label{HoldOrlicz} Let $\phi$ and $\phi^*$ be a pair of conjugate Young functions, and let $f \in L^\phi(\Omega)$ and $g \in L^{\phi^*}(\Omega)$. Then \[ \int_\Omega |f(x)g(x)|\, \dx \leq 2\, \|f\|_{L^\phi(\Omega)} \, \|g\|_{L^{\phi^*}(\Omega)}, \] where $\|\cdot\|_{L^\phi(\Omega)}$ and $\|\cdot\|_{L^{\phi^*}(\Omega)}$ denote the Luxemburg norms associated with $\phi$ and $\phi^*$, respectively. \end{proposition}
For further needs, we recall  the following well known  Moser Trudinger type imbedding theorem in the Orlicz Sobolev setting (\cite {Cianchi}, \cite{trudinger}).
\begin{theorem}\label{sobtru}
Let $\Omega$ be an open set in $\mathbb{R}^n$ with $|\Omega|<\infty$. Then, there exists a constant $\mathcal{C}_S=\mathcal{C}_S(n,|\Omega|)$ such that
    \begin{align}
        \label{SPnew1}
        ||w||_{{\rm EXP}_{\frac{n}{n-1}}(\Omega)} \leq \mathcal{C}_S\left( \int_\Omega |D w|^n \dx\right)^{\frac{1}{n}}
    \end{align}
    for every $w \in W^{1,n}_0(\Omega)$.
\end{theorem}
\iffalse
\tcb{We conclude this subsection recalling the following elementary inequality.
\begin{lemma}\label{elemnew}Let  $\Phi:[0,+\infty]\to [0,+\infty]$ be a strictly increasing function. Then the following inequality
$$st\le s\Phi^2\left(\frac{s}{\varepsilon}\right)+\varepsilon t\Phi^{-1}(\sqrt{t}),$$
for every $s,t\ge 0$ and every $\varepsilon>0$
\end{lemma}
\begin{proof}
If $t\le \Phi^2\left(\frac{s}{\varepsilon}\right)$, clearly we have
$$st\le s\Phi^2\left(\frac{s}{\varepsilon}\right).$$ Otherwise, if
$t> \Phi^2\left(\frac{s}{\varepsilon}\right)$
we have
$$\Phi^{-1}(\sqrt{t})>\frac{s}{\varepsilon}\,\,\Longleftrightarrow\,\, s<\varepsilon \Phi^{-1}(\sqrt{t})$$
and, in this case, we have
$$st<\varepsilon \Phi^{-1}(\sqrt{t}).$$
Joining the two cases we have the conclusion.
\end{proof}}
\fi
We now recall a classical embedding relating the $L^p$
 spaces to spaces of exponentially integrable functions.  
\begin{proposition}\label{normalpexp}
 Let $\Omega$ be an open set in $\mathbb{R}^n$ with $|\Omega|<\infty$. If $f\in \mathrm{EXP}_\alpha(\Omega)$  with $\alpha>0$, then there exists a constant
$C= {C(\alpha, |\Omega|)}>0$ such that
\[
\|f\|_{L^p(\Omega)}
\le C\, p^{1/\alpha}\, \|f\|_{ \mathrm{EXP}_\alpha(\Omega)},
\]
for every $p\ge 1$.
\end{proposition}
\begin{proof}
To simplify the notations, set
$A:=\|f\|_{ \mathrm{EXP}_\alpha(\Omega)},
$
so that, by the definition of the Luxemburg norm, 
\[
\int_{\Omega} \left(e^{(|f|/A)^\alpha}-1\right)\,\dx = 1.
\]
Using that $
e^t-1=\displaystyle{\sum_{k=1}^{\infty}\frac{t^k}{k!}},
$
with $t=(|f|/A)^\alpha$, we obtain
\[
\int_{{\Omega}} \sum_{k=1}^{\infty}
\frac{|f|^{\alpha k}}{A^{\alpha k}k!}\,\dx
= 1.
\]
It follows that 
\[
\frac{1}{A^{\alpha k}k!}
\int_{{\Omega}} |f|^{\alpha k}\,\dx
\le 1,
\]
for
every $k\ge1$, and so
\[
\|f\|_{L^{\alpha k}({\Omega})}^{\alpha k}
\le A^{\alpha k}k!,
\]
or equivalently,
\[
\|f\|_{L^{\alpha k}({\Omega})}
\le A\,(k!)^{1/(\alpha k)}.
\]
By Stirling's formula, there exists a constant $C>0$ such that
\[
(k!)^{1/k}\le Ck,
\]
for all $k\ge1$. Therefore,
\[
\|f\|_{L^{\alpha k}({\Omega})}
\le C^{1/\alpha}A\,k^{1/\alpha}.
\]
Setting $p=\alpha k$, we obtain
\[
\|f\|_{L^p({\Omega})}
\le C\,p^{1/\alpha}\,
\|f\|_{ \mathrm{EXP}_\alpha({\Omega})}
\]
for all $p$ of the form $p=\alpha k$.

The estimate for arbitrary $p\ge1$ follows from {by H\"older's inequality and the embedding of $L^q(\Omega)$ into $L^p(\Omega)$ for $p \le q$}
%$L^p$ norms on finite measure spaces (or, alternatively, by
%interpolation). 
Thus
\[
\|f\|_{L^p({\Omega})}
\le C\,p^{1/\alpha}\,
\|f\|_{{\rm EXP}\alpha({\Omega})},
\qquad p\ge1,
\]
where $C$ depends only on $\alpha$ {and $|\Omega|$.}
\end{proof}

{We conclude by recalling the  following
\begin{lemma}\label{GGPnew} 
%Let $\varphi:[0,+\infty)\to [0,+\infty)$ be  increasing and diverging at $\infty$. Then there exists $s_0>0$ such that 
%$$\frac{s}{\psi(s)}\le \varphi(s)\le 2\frac{s}{\psi(s)},\qquad \text{for every } s>s_0$$
%where for $\psi(t)=[(\varphi^{-1})^*]^{-1}(t)$.
Let $\Phi: [0, + \infty) \rightarrow [0, + \infty)$ an increasing function with $\Phi(0) = 0.$ Assume that
\[
t \mapsto \frac{\Phi(t)}{t}
\]
is positive and increasing on $[t_0, + \infty)$ with $t_0 > 0$ and that  verifies 
\[
\lim_{t \rightarrow 0} \frac{\Phi(t)}{t} = 0, \qquad \lim_{t \rightarrow + \infty} \frac{\Phi(t)}{t} = + \infty.
\]
Then
\[
\frac{s}{\Phi^{-1}(s)} \le \, (\Phi^*)^{-1}(s) \le \, 2 \frac{s}{\Phi^{-1}(s)} \qquad \forall s \in \, \Phi([t_0; \infty)), 
\]
where $\Phi^{*}$ is the Fenchel conjugate of $\Phi$ in the sense of \eqref{conjugate}. 
\end{lemma}
\color{black}
For the proof we refer to \cite[Lemma 2.4]{GGP}, see also \cite[Lemma 2.1]{GianGrePas}. {This lemma will be used at the end of the a priori estimate with the choice $\Phi = \mathcal{L}^{-1},$ being $\mathcal{L}$ the function appearing in assumption {(F6).}}

\iffalse
\colorbox{pink}{forse il lemma seguente non serve}
\begin{lemma}\label{elem}
For every $s,t\ge 0$ and every $\varepsilon>0$  the following inequality
$$st\le s\log\left(e+\frac{s}{\varepsilon}\right)+ \varepsilon \exp\left(2t\right)$$
holds.\end{lemma}}
\fi

\section{A Caccioppoli type inequality}
In this section we will prove a second order Caccioppoli inequality which underlies Moser’s iteration scheme.
The a priori assumption $u\in W^{1,\infty}_{\loc}(\Omega)\cap W^{2,\min\{2,p\}}_{\loc}(\Omega)$ allows us to  test the second variation of the functional $\mathcal{F}(u)$, with functions that are proportionals to the second derivatives of the minimizers.
\begin{lemma}
Let $u\in W^{1,\infty}_{\loc}(\Omega)\cap W^{2,\min\{2,p\}}_{\loc}(\Omega)$ be  a local minimizer of the functional $\mathcal{F}(u).$ Then the following second order Caccioppoli type inequality
\begin{eqnarray}\label{Caccioppoli}
 && \int_{\Omega} \eta^2 (1 + |Du|^2)^{\frac{p-2+\gamma}{2}}  |D^2 u|^2 \, \dx \cr\cr
&\le & C (1+\gamma^2)  \, \int_{\Omega} \eta^2  h^2(x)(1 + |Du|^2)^{\frac{p+\gamma}{2} } g^2((1 + |Du|^2)^{\frac{1}{2}}) \, \dx\cr\cr
&&+C \int_{\Omega} |D \eta|^2 \,  (1 + |Du|^2)^{\frac{p+\gamma}{2}} g((1 + |Du|^2)^{\frac{1}{2}})  \, \dx,
\end{eqnarray}
holds true for every {$\gamma \ge \, 0$} and for every $\eta\in C^1_0(\Omega).$
\end{lemma}
\begin{proof}
{ By the assumption 
\begin{equation}
\label{regE}
u \in  W^{1,\infty}_{\loc}(\Omega)\cap W^{2,\min\{2,p\}}_{\loc}(\Omega),
\end{equation}
we have that the following system
\begin{equation}
\label{second_variation}
\int_{\Omega} \left (\sum_{i,j=1}^n  F_{\xi_i \xi_j}(x, Du)u_{x_j x_s}D_{x_i} \varphi  \, dx + \sum_{i=1}^n F_{\xi_i x_s}(x, Du) D_{x_i} \varphi \right ) = 0,
\end{equation}
holds for all $s = 1, \dots, n$ and for all $\varphi \in W^{1,{\min\{2,p\}}}_0(\Omega)$. For $\eta \in \mathcal{C}^1_0(\Omega)$ and $\gamma \ge 0$, let us consider the function $\varphi=\eta^2
(1+|Du|^2)^{\frac{\gamma}{2}} u_{x_s}$,   so that
\begin{eqnarray*}
 D_{x_i} \varphi &=&2 \eta \eta_{x_i} (1 +|Du|^2)^{\frac{\gamma}{2}} u_{x_s} + \eta^2 (1 +|Du|^2)^{\frac{\gamma}{2}} u_{x_s x_i} \\
 && + \gamma \eta^2  (1 +|Du|^2)^{\frac{\gamma}{2}-1} |Du| D_{x_i}(|Du|) u_{x_s}.
\end{eqnarray*} 
Using $\varphi$ as test function in \eqref{second_variation}, we get
\begin{eqnarray*}
0 &=& 2\int_{\Omega} \eta\sum_{i,j=1}^n F_{\xi_i \xi_j}(x, Du) u_{x_j x_s}  \eta_{x_i} (1 +|Du|^2)^{\frac{\gamma}{2}} u_{x_s} \, \dx\\
&& +  \int_{\Omega} \eta^2\sum_{i,j=1}^n F_{\xi_i \xi_j}(x, Du) u_{x_j x_s}  (1 +|Du|^2)^{\frac{\gamma}{2}} u_{x_s x_i} \, \dx\\
&& + \gamma \int_{\Omega}\eta^2 \sum_{i,j=1}^n F_{\xi_i \xi_j}(x, Du) u_{x_j x_s}    (1 +|Du|^2)^{\frac{\gamma}{2}-1} |Du| D_{x_i}(|Du|) u_{x_s} \, \dx\\
&& +2 \int_{\Omega} \eta\sum_{i=1}^n F_{\xi_i x_s}(x, Du)  \eta_{x_i}  (1 +|Du|^2)^{\frac{\gamma}{2}} u_{x_s} \, \dx\\
&& +  \int_{\Omega} \eta^2\sum_{i=1}^n F_{\xi_i x_s }(x, Du)  (1 + |Du|^2)^{\frac{\gamma}{2}}  u_{x_s x_i} \, \dx\\
&& + \gamma \int_{\Omega} \eta^2\sum_{i=1}^n F_{\xi_i x_s}(x, Du)    (1 + |Du|^2)^{\frac{\gamma}{2}-1}  |Du| D_{x_i}(|Du|) u_{x_s} \, \dx .
\end{eqnarray*}
Summing all terms in the previous  equation   with respect to $s$ from 1 to $n$, and  denoting by ${I_1-I_6}$ the corresponding integrals, we obtain
\begin{equation}\label{start}
	I_2+I_3\le |{I}_1|+|{I}_4|+|{I}_5|+|{I}_6|.
\end{equation}
By  assumption (F1) and using that $D_{x_j}(|Du|)|Du|={\sum_{s=1}^n u_{x_j x_s} u_{x_s}}$, we can estimate the term $I_3$ as follows
\begin{eqnarray*}\label{I3}
{I}_3 &=& \gamma\int_{\Omega} \eta^2 \sum_{i,j,s =1}^n F_{\xi_i \xi_j}(x, Du) u_{x_j x_s} \left [   (1 + |Du|^2)^{\frac{\gamma}{2}-1}  D_{x_i}(|Du|)|Du| \right ] u_{x_s} \, \dx\cr\cr
%&\ge& 2 \gamma \, \int_{\Omega} \eta^2 |Du|^{2\gamma-1} \sum_{i,j,s =1}^n f_{\xi_i \xi_j}(x, Du) D_{x_i}(|Du|) u_{x_j x_s} u_{x_s} \, dx\\
%&=& 2 \gamma \, \int_{\Omega} \eta^2 |Du|^{2\gamma-1} \sum_{i,j =1}^n f_{\xi_i \xi_j}(x, Du) D_{x_i}(|Du|) \left(\sum_{s=1}^n u_{x_j x_s} u_{x_s} \right)\, dx\\
&\ge&  \gamma \, \int_{\Omega} \eta^2 (1+|Du|^2)^{\frac{\gamma}{2}-1}|Du| \sum_{i,j =1}^n F_{\xi_i \xi_j}(x, Du) D_{x_i}(|Du|) D_{x_j}(|Du|) \, \dx\cr\cr
&\ge&  \gamma \,  {\lambda} \int_{\Omega} \eta^2 (1+|Du|^2)^{\frac{\gamma}{2}-1}|Du|  |D(|Du|)|^2  (1 + |Du|^2)^{\frac{p-2}{2}}  \, \dx \ge 0.
\end{eqnarray*}
Hence,  \eqref{start} implies
\begin{eqnarray}\label{ristart}
	I_2\le |{I}_1|+|{I}_4|+|{I}_5|+|{I}_6|.
\end{eqnarray}
By  Cauchy-Schwarz and Young's inequalities, we have
\begin{eqnarray}\label{I1}
|{I}_1| &=& 2\left |\int_{\Omega}  \eta  (1 + |Du|^2)^{\frac{\gamma}{2}}   \sum_{i,j,s =1}^n F_{\xi_i \xi_j}(x, Du) u_{x_j x_s} \eta_{x_i} u_{x_s} \, \dx\right |\cr\cr
%& = &
%2\Bigg| \int_{\Omega}  \eta  (1 + |Du|^2)^{\gamma} \cr\cr
%&& \times \sum_{s=1}^n \left( \sum_{i,j=1}^n F_{\xi_i \xi_j}(x, Du) \eta_{x_i} \eta_{x_j} u_{x_s}^2 \right)^{1/2} \left(\sum_{i,j=1}^n
%F_{\xi_i \xi_j}(x, Du) u_{x_s x_i} \, u_{x_s x_j} \right)^{1/2} \, dx \Bigg| \cr\cr
&\le& 2\int_{\Omega}  \eta  (1 + |Du|^2)^{\frac{\gamma}{2}} \cr\cr
&& \times \left \{ \sum_{i,j,s =1}^n F_{\xi_i \xi_j}(x, Du) \eta_{x_i} \eta_{x_j} u_{x_s}^2\right \}^{1/2} \, \left \{ \sum_{i,j,s =1}^n F_{\xi_i \xi_j}(x, Du) u_{x_s x_i} \, u_{x_s x_j}\right \}^{1/2} \, \dx \cr\cr
&\le &  C   \int_{\Omega}  (1 + |Du|^2)^{\frac{\gamma}{2}}  \sum_{i,j,s =1}^n F_{\xi_i \xi_j}(x, Du) \eta_{x_i} \eta_{x_j} u_{x_s}^2 \, \dx\cr\cr
&&  + \frac{1}{2} \int_{\Omega} \eta^2  (1 + |Du|^2)^{\frac{\gamma}{2}}  \sum_{i,j,s =1}^n F_{\xi_i \xi_j}(x, Du) u_{x_s x_i} \, u_{x_s x_j}  \, \dx \cr\cr
%&\le &  C (\Lambda)  \int_{\Omega} \Lambda  (1 + |Du|^2)^{\frac{p-2}{2}+\gamma}  \sum_{i,j,s =1}^n \eta_{x_i} \eta_{x_j} u_{x_s}^2 \, dx \\
%&& + \frac{1}{4} \int_{\Omega} \eta^2  (1 + |Du|^2)^{\gamma}  \sum_{i,j,s =1}^n f_{\xi_i \xi_j}(x, Du) u_{x_s x_i} \, u_{x_s x_j} \, dx \\
%&\le &  C(n, \Lambda)  \int_{\Omega} |D\eta|^2  (1 + |Du|^2)^{\frac{p-2}{2}+\gamma}  \sum_{s =1}^n u_{x_s}^2 \, dx \\
%&& + \frac{1}{4} \int_{\Omega} \eta^2  (1 + |Du|^2)^{\gamma}  \sum_{i,j,s =1}^n f_{\xi_i \xi_j}(x, Du) u_{x_s x_i} \, u_{x_s x_j} \, dx \\
&\le& {C(\Lambda)} \, \int_{\Omega} |D \eta|^2 \,  (1 + |Du|^2)^{\frac{p+\gamma}{2}} {g\big((1+|Du|^2)^{\frac12}}\big) \, \dx \cr\cr
&& + \frac{1}{2} \int_{\Omega} \eta^2  (1 + |Du|^2)^{\frac{\gamma}{2}}  \, \sum_{i,j,s =1}^n F_{\xi_i \xi_j}(x, Du) u_{x_j x_s} u_{x_i x_s} \, \dx,
\end{eqnarray}
where, in the last inequality, we used assumption  (F2).
Using {once more} Young's inequality together with (F3) yields
\begin{eqnarray}\label{I4}
|{I}_4| &=& 2 \left | \int_{\Omega} \eta \, (1 + |Du|^2)^{\frac{\gamma}{2}}  \sum_{i,s =1}^n F_{\xi_i x_s}(x, Du) \eta_{x_i} u_{x_s} \, \dx \right |\cr\cr
%&{\le}&  2 \,  \int_{\Omega} \eta \, h(x) \,  (1 + |Du|^2)^{\frac{p-1+\gamma}{2}}g\big((1+|Du|^2)^{\frac12}\big)  \sum_{i,s=1}^n |\eta_{x_i} u_{x_s} | \, \dx \cr\cr
&\le&   C \,  \int_{\Omega} \eta \, |D\eta| |Du| \, h(x) \, (1 + |Du|^2)^{\frac{p-1+\gamma}{2}} g\big((1+|Du|^2)^{\frac12}\big)  \, \dx \cr\cr
&\le& C \int_{\Omega} \eta \, |D\eta| h(x) (1 + |Du|^2)^{\frac{p+\gamma}{2}}g\big((1+|Du|^2)^{\frac12}\big)\,\dx\cr\cr
&\le& C\int_{\Omega} |D\eta|^2  (1 + |Du|^2)^{\frac{p+\gamma}{2}}\,\dx\\\nonumber
&&+C\int_{\Omega} \eta^2 h^2(x) (1 + |Du|^2)^{\frac{p+\gamma}{2}}g^2\big((1+|Du|^2)^{\frac12}\big)\,\dx
\end{eqnarray}
and also
\begin{eqnarray}\label{I5}
|{I}_5| &=& \left| \int_{\Omega} \eta^2 (1 + |Du|^2)^{\frac{\gamma}{2}} \sum_{i,s=1}^n F_{\xi_i x_s }(x, Du)  u_{x_s x_i} \, \dx \right| \cr\cr
%&{\le}&  \int_{\Omega} \eta^2\,h(x) (1 + |Du|^2)^{\frac{p-1+\gamma}{2}}g\big((1+|Du|^2)^{\frac12}\big)   \left |\sum_{i,s=1}^n  u_{x_s x_i} \right | \, \dx \cr\cr
&\le&\,  \, \int_{\Omega} \eta^2 h(x) (1 + |Du|^2)^{\frac{p-1+\gamma}{2}}g\big((1+|Du|^2)^{\frac12}\big)  |D^2 u| \, \dx \cr\cr
&=& \int_{\Omega} \left[\eta^2 (1+|Du|^2)^{\frac{p-2+\gamma}{2}} |D^2 u|^2 \right]^{\frac{1}{2}} \left[\eta^2 (1+|Du|^2)^{\frac{ p+\gamma}{2}}g^2\big((1+|Du|^2)^{\frac12}\big)h^2(x)\right]^{\frac{1}{2}} \, \dx \cr\cr
&\le& \varepsilon \int_{\Omega} \eta^2 (1 + |Du|^2)^{\frac{p-2+\gamma}{2}}  |D^2 u|^2 \, \dx\\\nonumber&& + C_\varepsilon  \int_{\Omega} \eta^2 h^2(x)(1 + |Du|^2)^{\frac{ p+\gamma}{2}} g^2\big((1+|Du|^2)^{\frac12}\big) \, \dx,
\end{eqnarray}
where $\varepsilon>0$ will be chosen later.
Arguing similarly, we get
\begin{eqnarray}\label{I6}
|{I}_6| &=&  \left | \gamma \, \int_{\Omega} \sum_{i,s=1}^n F_{\xi_i x_s}(x, Du) \eta^2  (1 + |Du|^2)^{\frac{\gamma}{2}-1}  |Du| D_{x_i}(|Du|) u_{x_s} \, \dx \right |\cr\cr
%&=&2 \gamma \, \int_{\Omega} \eta^2 (1 + |Du|^2)^{\gamma-1}  |Du| \sum_{i,s=1}^n f_{\xi_i x_s}(x, Du)  D_{x_i}(|Du|) u_{x_s} \, dx\\
&\le&  \, \gamma \, \int_{\Omega} \eta^2 (1+|Du|^2)^{\frac{\gamma-1}{2}}  {|Du|}\left |\sum_{i,s=1}^n F_{\xi_i x_s}(x, Du)  D_{x_i}(|Du|) u_{x_s} \, \dx \right |\cr\cr
%&\le&  \, \gamma \, \int_{\Omega} \eta^2  h(x)(1 + |Du|^2)^{\frac{p-2+\gamma}{2}} \,  g\big((1+|Du|^2)^{\frac12}\big)\,\left | \sum_{i,s=1}^n D_{x_i}(|Du|) u_{x_s} \, dx \right |\cr\cr
%&\le& 2 \, n \, \gamma \, \Lambda \, \int_{\Omega} \eta^2 (1 + |Du|^2)^{\gamma-\frac{1}{2}+\frac{q-1}{2}}  |D(|Du|)| |Du| \, dx\\
&\le& C \, \gamma \, \int_{\Omega} \, \eta^2 (1 + |Du|^2)^{\frac{p-1+\gamma}{2}}  g\big((1+|Du|^2)^{\frac12}\big)|D^2 u| \, h(x)\, \dx\cr\cr
%&=& C \gamma\int_{\Omega} \left( \eta  h(x) (1+|Du|^2)^{\frac{2q-p}{4}+\frac{\gamma}{2}} \right)  \left( \eta |D^2u| (1+|Du|^2)^{\frac{p-2}{4}+\frac{\gamma}{2}} \right) \, dx\cr\cr
&\le& \varepsilon \int_{\Omega} \eta^2 |D^2u|^2 (1 + |Du|^2)^{\frac{p-2+\gamma}{2}}   \, \dx\\\nonumber
&&+ C_\varepsilon \gamma^2 \,  \, \int_{\Omega} \eta^2  h^2(x)(1 + |Du|^2)^{\frac{p+\gamma}{2} }g^2\big((1+|Du|^2)^{\frac12}\big)  \, \dx,
\end{eqnarray}
where  the constants $C$ and $C_{\varepsilon}$ are independent of $\gamma$.}
Inserting  \eqref{I1}, \eqref{I4}, \eqref{I5} and \eqref{I6} in \eqref{ristart} we obtain
\begin{eqnarray*}
&& \int_{\Omega} \eta^2 (1 +|Du|^2)^{\frac{\gamma}{2}}{\sum_{i,j,s=1}^n} F_{\xi_i \xi_j}(x, Du) u_{x_j x_s}  u_{x_s x_i} \, \dx\cr\cr
&\le&  \frac{1}{2} \int_{\Omega} \eta^2  (1 + |Du|^2)^{\frac{\gamma}{2}}  \, \sum_{i,j,s =1}^n F_{\xi_i \xi_j}(x, Du) u_{x_j x_s} u_{x_i x_s} \, \dx\cr\cr
&&+2\varepsilon \int_{\Omega} \eta^2 (1 + |Du|^2)^{\frac{p-2+\gamma}{2}}  |D^2 u|^2 \, \dx \cr\cr
&&+ C_\varepsilon (1+\gamma^2)  \, \int_{\Omega} \eta^2  h^2(x)(1 + |Du|^2)^{\frac{p+\gamma}{2} } g^2\big((1+|Du|^2)^{\frac12}\big) \, \dx\cr\cr
&&+ C(\varepsilon, \Lambda)\int_{\Omega} |D \eta|^2 \,  (1 + |Du|^2)^{\frac{p+\gamma}{2}} {g\big((1+|Du|^2)^{\frac12}\big)}  \, \dx,
\end{eqnarray*}
where, in the last line we used that $g(1+t)\ge 1.$
Reabsorbing the first integral in the right hand side by the left hand side, we get
\begin{eqnarray*}
&& \frac{1}{2}\int_{\Omega} \eta^2 (1 +|Du|^2)^{\frac{\gamma}{2}}{\sum_{i,j,s=1}^n} F_{\xi_i \xi_j}(x, Du) u_{x_j x_s}  u_{x_s x_i} \, \dx\cr\cr
 &\le &2\varepsilon \int_{\Omega} \eta^2 (1 + |Du|^2)^{\frac{p-2+\gamma}{2}}  |D^2 u|^2 \, \dx \cr\cr
&&+ C_\varepsilon (1+\gamma^2) \,  \int_{\Omega} \eta^2  h^2(x)(1 + |Du|^2)^{\frac{p+\gamma}{2}} g^2\big((1+|Du|^2)^{\frac12}\big) \, \dx\cr\cr
&&+ C(\varepsilon, \Lambda) \int_{\Omega} |D \eta|^2 \,  (1 + |Du|^2)^{\frac{p+\gamma}{2}} {g\big((1+|Du|^2)^{\frac12}\big)}  \, \dx.
\end{eqnarray*}
Using assumption (F1) in the left hand side of previous estimate, we obtain
\begin{eqnarray*}
&& \frac{{\lambda}}{2}\int_{\Omega} \eta^2 (1 + |Du|^2)^{\frac{p-2+\gamma}{2}}  |D^2 u|^2 \, \dx \cr\cr
 &\le &2\varepsilon \int_{\Omega} \eta^2 (1 + |Du|^2)^{\frac{p-2+\gamma}{2}}  |D^2 u|^2 \, \dx \cr\cr
&&+ C_\varepsilon (1+\gamma^2) \, \int_{\Omega} \eta^2  h^2(x)(1 + |Du|^2)^{\frac{p+\gamma}{2} }g^2\big((1+|Du|^2)^{\frac12}\big)  \, \dx\cr\cr
&&+C_\varepsilon \int_{\Omega} |D \eta|^2 \,  (1 + |Du|^2)^{\frac{p+\gamma}{2}} {g\big((1+|Du|^2)^{\frac12}\big)} \, \dx.
\end{eqnarray*}
Choosing {for instance} $\varepsilon=\frac{{\lambda}}{8}$, we can reabsorb the first integral in the right hand side by the left hand side thus getting
\begin{eqnarray*}
&& \int_{\Omega} \eta^2 (1 + |Du|^2)^{\frac{p-2+\gamma}{2}}  |D^2 u|^2 \, \dx \cr\cr
 &\le & C (1+\gamma^2)  \, \int_{\Omega} \eta^2  h^2(x)(1 + |Du|^2)^{\frac{p+\gamma}{2} }  g^2\big((1+|Du|^2)^{\frac12}\big)\, \dx\cr\cr
&&+C \int_{\Omega} |D \eta|^2 \,  (1 + |Du|^2)^{\frac{p+\gamma}{2}} {g\big((1+|Du|^2)^{\frac12}}\big) \, \dx,
\end{eqnarray*}
with a constant $C=C(\lambda,\Lambda)>0$ and this concludes the proof.
\end{proof}

\medskip

\section{The a priori estimate {in the case $\mathbf{n\ge 3}$}}

This section is devoted to the proof of a uniform a priori estimate for the Lipschitz norm of the minimizers of the functional at \eqref{functional} {which relies on a Moser type iteration}. Here we face  the case $n\ge 3$, while the case $n=2$ will be treated in the subsequent section, since it requires the use of a different Sobolev imbedding Theorem. The precise statement is the following
\begin{theorem}\label{thmapriori}
Let $u\in W^{1,\infty}_{\loc}(\Omega)\cap W^{2,{\min\{2,p\}}}_{\loc}(\Omega)$ be  a local minimizer of the functional $\mathcal{F}(u,\Omega)$ {at \eqref{functional}}, under the assumptions {{{\rm (F1)--(F6)}}.} Then the following estimate \begin{eqnarray}\label{apriori}
	||V(Du)||_{L^{\infty}(B_{\rho_0})} &\le&   \bigg(\mathcal{L}^{-1}\bigg)^{**}\Bigg(C\left(\int_{B_{  R_0}}V(Du)^{p}\, \dx\right)^{\frac{1}{p}}\Bigg)\cr\cr
&&+\frac{C}{(R_0- \rho_0 )^{{\frac{4 \, 2^*}{p (2^* - 2)}}}}\left(\int_{B_{  R_0}}V(Du)^{p}
\, \dx\right)^{\frac{1}{p}}
%\cr\cr&& {+ \Big(\tilde g^{-1}\Big)^*\Bigg(\frac{C}{(\bar R-\baIr \rho )^{\frac{n}{p}}}\left(\int_{B_{\bar R}}V(Du)^{p}
%\, dx\right)^{\frac{1}{p}}\Bigg)}
\end{eqnarray}
holds for all balls $B_{\rho_0}\subset B_{R_0}\Subset\Omega$, {with constants $C$ depend on $n,p,\lambda, \Lambda, h$ in the case $n \ge 3$ and on $p,\lambda, \Lambda, h, 2^*, \theta, \mu, c_{\Phi}, |B_{R_0}|$ in the case $n = 2$} and where $\mathcal{L}(t)$ is the function appearing in {{\rm(F6)}}.
%$\tilde{\Psi}(t)=\Psi^{\frac{n}{2p}}(t)$ \tcv{and $\tilde g(t)= g^{\frac{n}{2p}}(t)$.}
\end{theorem}
\begin{proof}[Proof in the case $\mathbf{n\ge 3}$]
Let us fix a ball $B_{R_0}\Subset \Omega$ and radii $\rho_{0}<\bar\rho<\rho<R<\bar R<R_{0}.$
 We choose in  \eqref{Caccioppoli} a cut-off function $\eta\in C^1_0(B_{R})$ such that $\eta=1$ on $B_{\rho}$ and $|D\eta|\le \frac{C}{R-\rho}$, {thus getting}
\begin{eqnarray}\label{estprimadellasob}
&& \int_{B_{R}} \eta^2 (1 + |Du|^2)^{\frac{p-2 + \gamma}{2}}  |D^2 u|^2 \, \dx \cr\cr
&\le & C (1+\gamma^2)  \, \int_{B_{R}}\eta^2  h^2 (1 + |Du|^2)^{\frac{p + \gamma}{2}} g^2\big({(1+|Du|^2)^{\frac12}}\big) \, \dx\cr\cr
&&\quad+  \frac{C}{(R-\rho)^2}  \, \int_{B_{R}}   (1 + |Du|^2)^{\frac{p + \gamma}{2}}  {g\big({(1+|Du|^2)^{\frac12}}\big)} \, \dx.\label{primastima}
\end{eqnarray}
The Sobolev  inequality yields
\begin{eqnarray*}
&& \left ( \int_{B_{R}}  \eta^{2^*}(1 + |Du|^2)^{(\frac{p + \gamma}{4}) 2^*}\, \dx\right )^{\frac{2}{2^*}}\\
& \le& \, C \, \int_{B_{R}} \Big|D (\eta(1 + |Du|^2)^{\frac{p + \gamma}{4}})\Big|^2 \, \dx 	\\
&\le&C(1+\gamma^2)\int_{B_{R}} \eta^2 (1 + |Du|^2)^{\frac{p-2 + \gamma}{2}}  |D^2 u|^2\,\dx+C \int_{B_{R}} |D \eta|^2 \,  (1 + |Du|^2)^{\frac{p + \gamma}{2}}  \, \dx,
\end{eqnarray*}
where $2^*$  denotes, as usual, the Sobolev conjugate exponent of $2$, i.e., since $n>2$, $2^*=\frac{2n}{n-2}$.
%$$ 2^*=\begin{cases}
	%\displaystyle{\frac{2n}{n-2}}\qquad\qquad\qquad\qquad\qquad \text{if}\,\,n\ge3\cr\cr
	%\text{any finite exponent}\, \qquad\quad\qquad \text{if}\,\,n=2.
	%\end{cases}$$
	Using estimate \eqref{primastima} to control the first integral in the right hand side of the previous inequality {and exploiting once more the fact that $g(1 + t) \ge \, 1,$} we obtain
\begin{eqnarray}\label{stimapreite}
&&\left ( \int_{B_{R}}  \eta^{2^*}(1 + |Du|^2)^{(\frac{p + \gamma}{2}) \frac{2^*}{2}}\, \dx\right )^{\frac{2}{2^*}} 	\cr\cr
&\le& C {(1+\gamma^2)^2} \int_{B_{R}}\eta^2 h^2 (1 + |Du|^2)^{\frac{p + \gamma}{2}} g^2\big({(1+|Du|^2)^{\frac12}}\big) \, \dx \cr\cr
&&+\frac{C(1+\gamma^2)}{(R - \rho)^{2}} \int_{B_{R}}  \,  (1 + |Du|^2)^{\frac{p + \gamma}{2}} {g\big({(1+|Du|^2)^{\frac12}}\big)}  \, \dx.
\end{eqnarray}	
{From now on, to shorten the notations, we set} 
\begin{equation}
\label{V(Du)}
V(Du)=(1+|Du|^2)^{\frac{1}{2}},
\end{equation}{and we write} \eqref{stimapreite} as follows
\begin{eqnarray}\label{stimapreitebis}
&&\left ( \int_{B_{R}}  \eta^{2^*}V(Du)^{(p+\gamma)\frac{2^*}{2}}\, \dx\right )^{\frac{2}{2^*}} 	\cr\cr
&\le& C {(1+\gamma^2)^2}  \int_{B_R} \eta^2 h^2 V(Du)^{p + \gamma} g^2({V(Du)})  \, \dx\cr\cr
&&\quad+\frac{C(1+\gamma^2)}{(R - \rho)^{2}} \int_{B_{R}}  \,  V(Du)^{p + \gamma}  {g({V(Du)})} \, \dx \cr\cr
&\le &  \varepsilon^2 {(1+\gamma^2)^2}  \int_{{B_R \cap} \{h < \varepsilon\}}\eta^2 V(Du)^{p + \gamma}  g^2({V(Du)}) \, \dx\cr\cr
&&+ C {(1+\gamma^2)^2}  \int_{{B_R \cap} \{h \ge\varepsilon\}}\eta^2
\frac{h^2}{\Phi^2(h)} \Phi^2(h) V(Du)^{p + \gamma} g^2({V(Du)}) \, \dx \cr \cr 
&&+\frac{C (1+\gamma^2)}{(R - \rho)^{2}} \int_{B_{R}}  \,  V(Du)^{p + \gamma}{g({V(Du)})}  \, \dx \cr \cr 
&\le& \varepsilon^2 {(1+\gamma^2)^2}  g^2(\|{V(Du)}\|_{{L^\infty(B_{R_0})}}) \int_{B_R} V(Du)^{p + \gamma} \, \dx\cr\cr && +\frac{C(1+\gamma^2) \, {g(\|{V(Du)}\|_{{L^\infty(B_{R_0})}})} }{(R - \rho)^{2}} \int_{B_{R}}  \,  V(Du)^{p + \gamma}  \, \dx\cr \cr 
&& + C\frac{{(1+\gamma^2)^2} }{\Phi^2(\varepsilon)} g^2\big(\|{V(Du)}\|_{{L^\infty(B_{R_0})}}\big) \int_{B_R}\eta^2 h^2 {\Phi^2(h)} V(Du)^{p + \gamma} \,\dx\cr\cr
&=:& {\mathrm I} + \mathrm{II} + \mathrm{III},
%&\le& C\frac{(1+\gamma^4)}{(R - \rho)^{2}}||V(Du)||^{2(q-p)}_{L^\infty(B_R)}\left [\int_{B_R}V(Du)^{(p+2\gamma)m}  \, dx \right ]^{\frac{1}{m}}
\end{eqnarray}
where we used the monotonicity of the functions $g$ and  $\Phi$, the a priori assumption $Du\in L^\infty(B_{R_0})$, and $\varepsilon>0$ is a parameter that will be chosen later.
We estimate the integral $III$ by H\"older's inequality as follows
\begin{eqnarray*}
\mathrm{III} &\le& C\frac{{(1+\gamma^2)^2} }{\Phi^2(\varepsilon)} g^2\big(\|{V(Du)}\|_{{L^\infty(B_{R_0})}}\big) \left (\int_{B_R}h^n \Phi(h)^{n} \, \dx \right )^{2/n} \cr\cr
&&\qquad\qquad\cdot\left (\int_{B_R}\eta^{2^*}V(Du)^{\frac{(p + \gamma) n}{n-2}} \, \dx \right )^{\frac{n-2}{n}},
\end{eqnarray*}
which is legitimate by virtue of assumption {(F4),} that holds true with $\theta = n$ since we are in the case $n \ge \, 3.$
Setting
\begin{equation}
\label{def_Hn}
{H_n := \left (\int_{B_{R_0}} h^n \Phi(h)^{n} \, \dx \right )^{1/n},}
\end{equation}
the previous estimate reduces to
\begin{eqnarray*}
\mathrm{III} &\le& C\frac{{H_n^2 (1+\gamma^2)^2} }{\Phi^2(\varepsilon)} g^2\big(\|{V(Du)}\|_{{L^\infty(B_{R_0})}} \big)\left (\int_{B_R}\eta^{2^*} V(Du)^{\frac{(p + \gamma) n}{n-2}} \, \dx \right )^{\frac{n-2}{n}}.
\end{eqnarray*}
At this point we choose $\varepsilon$ so that
\[
C\frac{{H_n^2 \, (1+\gamma^2)^2} }{\Phi^2(\varepsilon)} g^2\big(\|{V(Du)}\|_{{L^\infty(B_{R_0})}}\big) = \frac{1}{2}.
\]
This yields
\[
\Phi^2(\varepsilon) = 2C {H_n^2 (1+\gamma^2)^2}  g^2\big(\|{V(Du)}\|_{{L^\infty(B_{R_0})}}\big) \quad \Longleftrightarrow \quad \Phi(\varepsilon) = \tilde C {H_n}(1+\gamma^2) g\big(\|{V(Du)}\|_{{L^\infty(B_{R_0})}}\big)
\]
$$\quad \Longrightarrow \quad \varepsilon = \Phi^{-1} \Big({\tilde{C}} {H_n}(1+\gamma^2) g\big(\|{V(Du)}\|_{{L^\infty(B_{R_0})}}\big)\Big)$$
and, with this choice, {shortening $g\big(\|V(Du)\|_{{L^\infty(B_{R_0})}}\big)$ with $g\big(\|V(Du)\|_{\infty}\big)$}, estimate \eqref{stimapreitebis} becomes
\begin{eqnarray*}
&&\left (\int_{B_{{R}}} \eta^{2^*}V(Du)^{\frac{(p + \gamma) n}{n-2}} \,\dx \right )^{\frac{n-2}{n}}\cr\cr
&\le& \frac{1}{2} \left (\int_{B_R}\eta^{2^*} V(Du)^{\frac{(p + \gamma) n}{n-2}} \, \dx \right )^{\frac{n-2}{n}}+\frac{C \, {(1+\gamma^2)g(\|{V(Du)}\|_{\infty})}}{(R - \rho)^{2}} \int_{B_{R}}  \,  V(Du)^{p + \gamma}  \, \dx\cr\cr
&&\quad +C\left[(1+\gamma^2)g(\|{V(Du)}\|_{\infty}) \Phi^{-1} \Big( \tilde C {H_n}(1+\gamma^2) g(\|{V(Du)}\|)_{\infty}\Big)\right]^2 \int_{B_R} V(Du)^{p + \gamma} \, \dx.
\end{eqnarray*}
Reabsorbing the first integral in the right hand side by the left hand side, and recalling that $\eta=1$ on $B_\rho$, we deduce that
\begin{eqnarray}\label{baseiter}
&&\left (\int_{B_{\rho}} V(Du)^{\frac{(p + \gamma) n}{n-2}} \,\dx \right )^{\frac{n-2}{n}} \cr\cr
&\le& C\left[(1+\gamma^2)g\big(\|{V(Du)}\|_{\infty} \big)\Phi^{-1} \Big(\tilde C {H_n} (1+\gamma^2) g(\|{V(Du)}\|_{\infty})\Big)\right]^2 \int_{B_R} V(Du)^{p + \gamma} \, \dx\cr\cr
&&\qquad+C\frac{ (1+\gamma^2)\, {g(\|{V(Du)}\|_{\infty})} }{(R - \rho)^{2}} \int_{B_{R}}  \,  V(Du)^{p + \gamma}  \, \dx\cr\cr
&\le& 
{C(1+\gamma^4)g^2\big(\|{V(Du)}\|_{\infty} \big){\Phi^{-1} \Big({\tilde C^2} {H_n^2}(1+\gamma^4)\Big)\Phi^{-1} \Big( g^2(\|{V(Du)}\|_{\infty})\Big)} \int_{B_R} V(Du)^{p + \gamma} \, \dx}\cr\cr
&&\qquad+C\frac{ \, {(1+\gamma^2)g(\|{V(Du)}\|_{\infty})}}{(R - \rho)^{2}} \int_{B_{R}}  \,  V(Du)^{p + \gamma}  \, \dx.
\end{eqnarray}
where, in the second line of estimate, we used Lemma \ref{elem1} that is legitimate thanks to $\mathrm{(F5)}_1$. %and where, without loss of generality, we supposed $R-\rho\le 1$.
Note that the previous estimate holds true for every $\rho<R$ such that $\bar\rho\le  \rho< R\le \bar R$ and,  by virtue of the a priori assumption $u\in W^{1,\infty}_{\mathrm{loc}}(\Omega)$, for every $\gamma\ge 0$. Hence, defining the decreasing sequence of radii by setting
$$\rho_i=\bar \rho+\frac{ \bar R- \bar\rho}{2^i}$$
and the following increasing sequence of exponents
$$p_0=p \qquad\qquad {p_{i}}={p_{i-1}}\frac{2^*}{2}=p_0\left(\frac{2^*}{2}\right)^{i},$$
 we may iterate estimate \eqref{baseiter}  on the concentric balls $B_{\rho_i}$  with exponents $p_i$,
thus obtaining
\begin{eqnarray}\label{stimapreitepassoi}
&&\left ( \int_{B_{\rho_{i+1}}}  V(Du)^{p_{i+1}}\, \dx\right )^{\frac{1}{p_{i+1}}} 	\cr\cr
&\le& \displaystyle{\prod_{j=0}^{i}}\left[C\,p_j^{4}\,g^2(\|V(Du)\|_{\infty})\Phi^{-1} \big({\tilde C^2} {H_n^2}\,p_j^4\big)\Phi^{-1} \Big(g^2\big(\|V(Du)\|_{\infty}\big)\Big)+\frac{Cp_j^{2}\,g(\|V(Du)\|_{\infty})}{(\rho_j - \rho_{j+1})^{2}}\right]^{\frac{1}{p_j}}
\cr\cr
&&\qquad\qquad\cdot\left(\int_{B_{\bar R}}V(Du)^{p}  \, \dx\right)^{\frac{1}{p}}\cr\cr
&=& \displaystyle{\prod_{j=0}^{i}}\left[C\,p_j^{4}\,g^2(\|V(Du)\|_{\infty})\Phi^{-1} \big({\tilde C^2} {H_n^2}\,p_j^4\big)\Phi^{-1} \Big(g^2\big(\|V(Du)\|_{\infty}\big)\Big)+\frac{C 4^{j+1}p_j^{2}\,g(\|V(Du)\|_{\infty})}{(\bar R - \bar\rho)^{2}}\right]^{\frac{1}{p_j}}
\cr\cr
&&\qquad\qquad\cdot\left(\int_{B_{\bar R}}V(Du)^{p}  \, \dx\right)^{\frac{1}{p}}.
%\cr\cr
%&=& \displaystyle{\prod_{j=0}^{i}}\left({\frac{C4^{j+1} p_j^{4}{\left(\Phi^{-1} \Big(\tilde C H^2p_j^4g^2(\|Du\|_{\infty})\Big)\right)^2}}{(\bar R-\bar \rho )^{2}}} \right)^{\frac{1}{p_j}}\displaystyle{\prod_{j=0}^{i}}\left(g^2(\|Du\|_{\infty}) {+ g(\|Du\|_{\infty})}  \right)^{\frac{1}{p_j}}\cr\cr
%&&\qquad\qquad\cdot\left(\int_{B_{\rho_0}}V(Du)^{p}  \, dx\right)^{\frac{1}{p}}.
\end{eqnarray}
Since
\begin{eqnarray*}
\Theta_i &=& \displaystyle{\prod_{j=0}^{i}}\left[C\,p_j^{4}\,g^2(\|V(Du)\|_{\infty})\Phi^{-1} \big({\tilde C^2}{H_n^2}\,p_j^4\big)\Phi^{-1} \Big(g^2\big(\|V(Du)\|_{\infty}\big)\Big)+\frac{C 4^{j+1}p_j^{2}\,g(\|V(Du)\|_{\infty})}{(\bar R - \bar\rho)^{2}}\right]^{\frac{1}{p_j}}\\
&=&\exp\Bigg(\sum_{j=0}^i \frac{1}{p_j}\log\Bigg(C\,p_j^{4}\,g^2(\|V(Du)\|_{\infty})\Phi^{-1} \big({\tilde C^2} {H_n^2}\,p_j^4\big)\Phi^{-1} \Big(g^2\big(\|V(Du)\|_{\infty}\big)\Big)\cr\cr
&&\qquad\qquad\qquad+\frac{C 4^{j+1}p_j^{2}\,g(\|V(Du)\|_{\infty})}{(\bar R - \bar\rho)^{2}}\Bigg)\Bigg)\\
&\le &\exp\Bigg(\sum_{j=0}^i \frac{1}{p_j}\log\Bigg(C\Big(p_j^{4}\Phi^{-1} \big({\tilde C^2} {H_n^2}\,p_j^4\big)+4^{j+1}p_j^{2}\Big)\bigg(\,g^2(\|V(Du)\|_{\infty})\Phi^{-1} \Big(g^2\big(\|V(Du)\|_{\infty}\big)\Big)\cr\cr
&&\qquad\qquad\qquad+\frac{g(\|V(Du)\|_{\infty})}{(\bar R - \bar\rho)^{2}}\bigg)\Bigg)\Bigg)\\
&=& \exp\Bigg(\sum_{j=0}^i \frac{1}{p_j}\log\Bigg(C\Big(p_j^{4}\Phi^{-1} \big({\tilde C^2}{H_n^2}p_j^4\big)+4^{j+1}p_j^{2}\Big)\Bigg)\\
&& \cdot\exp\left(\sum_{j=0}^{i} \frac{1}{p_j}\log \Bigg(\,g^2(\|V(Du)\|_{\infty})\Phi^{-1} \Big(g^2\big(\|V(Du)\|_{\infty}\big) {\Big )} +\frac{g(\|V(Du)\|_{\infty})}{(\bar R - \bar\rho)^{2}}\Bigg)\right)\\
&\le & \exp\Bigg(\sum_{j=0}^i \frac{1}{p_j}\log\Big(1+C\,p_j^{4}\Phi^{-1} \big({\tilde C^2} {H_n^2}p_j^4\big)\Big)\Bigg)\cdot\exp\Bigg(\sum_{j=0}^i \frac{1}{p_j}\log\Big(1+C\,4^{j+1}p_j^{2}\Big)\Bigg)\\
&& \cdot\exp\left(\sum_{j=0}^{i} \frac{1}{p_j}\log \Bigg(\,g^2(\|V(Du)\|_{\infty})\Phi^{-1} \Big(g^2\big(\|V(Du)\|_{\infty}\big) {\Big )}+\frac{g(\|V(Du)\|_{\infty})}{(\bar R - \bar\rho)^{2}}\Bigg)\right)\\
&\le & \exp\Bigg(\sum_{j=0}^{+\infty} \frac{1}{p_j}\log\Big(1+C\,p_j^{4}\Big)\Bigg)\cdot\exp\Bigg(\sum_{j=0}^{+\infty} \frac{1}{p_j}\log\Big(1+\Phi^{-1} \big({\tilde C^2} {H_n^2}p_j^4\big)\Big)\Bigg)\\
&&\cdot\exp\Bigg(\sum_{j=0}^{+\infty} \frac{1}{p_j}\log\Big(1+4^{j+1}\Big)\Bigg)\cdot\exp\Bigg(\sum_{j=0}^{+\infty} \frac{1}{p_j}\log\Big(1+Cp_j^{2}\Big)\Bigg)\\
&& \cdot\exp\left(\sum_{j=0}^{+\infty} \frac{1}{p_j}\log\Bigg(\,g^2(\|V(Du)\|_{\infty})\Phi^{-1} \Big(g^2\big(\|V(Du)\|_{\infty}\big)\Big)+\frac{g(\|V(Du)\|_{\infty})}{(\bar R - \bar\rho)^{2}}\Bigg)\right)\\
&=:&\Theta_1\cdot \Theta_2\cdot \Theta_3\cdot \Theta_4\cdot \Theta_5,
\end{eqnarray*}
where we used that $\log(a+b)\le \log (1+a)+\log (1+b)$ {and also that $\log(1 + ab) \le \, \log(1 + a) \, \log(1 + b)$} for $a,b>0.$ 
Note that by the definition of $p_j$ and since $\sum_{j=0}^{{\infty}}\frac{1}{p_j}=\frac{2^*}{p(2^*-2)}= \frac{n}{ 2p} $ 
we have
$$\Theta_5\le \Bigg(\,g^2(\|V(Du)\|_{\infty})\Phi^{-1}\Big(g^2\big(\|V(Du)\|_{\infty}\big) {\Big )} +\frac{g(\|V(Du)\|_{\infty})}{(\bar R - \bar\rho)^{2}}\Bigg)^{\frac{n}{2p}}$$
and also, that
$$\Theta_1\cdot  \Theta_3\cdot \Theta_4\le c(n,p).$$
Moreover  assumption {(F5)$_2$} yields that
$$ \Theta_2\le \exp\left(c_\Phi\sum_{j=0}^{+\infty} \frac{1}{p_j}\Big({\tilde C^2} {H_n^2} p_j^4\Big)^\mu\right)\le c(n,p, {H_n},\lambda,\Lambda,\mu,c_\Phi),$$
since $\mu<\frac14$ {(this estimates holds when ${\tilde C^2} {H_n^2} p_j^4 > 1$ which is not restrictive to assume).}
Combining all the above estimates, we conclude that
\begin{eqnarray*}
 \Theta_i&\le&  C\Bigg(\,g^2(\|V(Du)\|_{\infty})\Phi^{-1}\Big(g^2\big(\|V(Du)\|_{\infty}\big)+\frac{g(\|V(Du)\|_{\infty})}{(\bar R - \bar\rho)^{2}}\Bigg)^{\frac{n}{2p}}\cr\cr
&\le& C\Bigg(\,g^{\frac{n}{p}}(\|V(Du)\|_{\infty})\Big[\Phi^{-1}\Big(g^2\big(\|V(Du)\|_{\infty}\big)\Big]^{\frac{n}{2p}}+\frac{g^{\frac{n}{2p}}(\|V(Du)\|_{\infty})}{(\bar R - \bar\rho)^{\frac{n}{p}}}\Bigg)
\end{eqnarray*}
for every $i\in \mathbb{N}$, with $C=C(n,p, {H_n},\lambda,\Lambda, \mu, c_\Phi)$.
%\tcb{To justify last estimate, we used that
%$$(1+a)(1+a\Phi^{-1}(a^2))\le \begin{cases}
   % 2(1+\Phi^{-1}(1))=:C\qquad\qquad\qquad\qquad \text{if}\,\, a\le 1\cr\cr
    %2a(1+a\Phi^{-1}(a^2))=2a+2a^2\Phi^{-1}(a^2)\quad \text{if}\,\, a> 1
%\end{cases} $$ that yields
%$$ (1+a)(1+a\Phi^{-1}(a^2))\le C(1+a+a^2\Phi^{-1}(a^2))$$}
Therefore, we can let $i\to \infty$ in \eqref{stimapreitepassoi} thus getting
\begin{eqnarray}\label{stimaquasifin}
	||V(Du)||_{L^{\infty}(B_{\bar\rho})} &\le&  C\Bigg(\,g^{\frac{n}{p}}(\|V(Du)\|_{\infty})\Big[\Phi^{-1}\Big(g^2\big(\|V(Du)\|_{\infty}\big)\Big]^{\frac{n}{2p}}+\frac{g^{\frac{n}{2p}}(\|V(Du)\|_{\infty})}{(\bar R - \bar\rho)^{\frac{n}{p}}}\Bigg)\cr\cr
&&\qquad\qquad\cdot \left(\int_{B_{ \bar R}}V(Du)^{p}
\, \dx\right)^{\frac{1}{p}}\cr\cr
 &=& C\left (\Psi^{\frac{n}{2p}}\big(\|V(Du)\|_{\infty}\big) + \frac{ g^{\frac{n}{2p}}(\|V(Du)\|_{\infty})}{(\bar R- \bar\rho )^{\frac{n}{p}}} \right )
 \cdot \left(\int_{B_{ \bar R}}V(Du)^{p}
\, \dx\right)^{\frac{1}{p}}\cr\cr
 &\le & C \left (\Psi^{\frac{n}{2p}}\big(\|V(Du)\|_{\infty}\big) + g^{\frac{n}{p}}(\|V(Du)\|_{\infty})+ \frac{{1}}{(\bar R- \bar\rho )^{\frac{2n}{p}}} \right )
  \cr\cr 
&&\qquad\qquad\cdot \left(\int_{B_{ \bar R}}V(Du)^{p}
\, \dx\right)^{\frac{1}{p}}\cr\cr 
&= & C\left (\Psi^{\frac{2^*}{p(2^*-2)}}\big(\|V(Du)\|_{\infty}\big) + g^{\frac{2\cdot2^*}{p(2^*-2)}}(\|V(Du)\|_{\infty})+ {\frac{{1}}{(\bar R- \bar\rho )^{\frac{2n}{p}}}}\right )
  \cr\cr 
&&\qquad\qquad\cdot \left(\int_{B_{ \bar R}}V(Du)^{p}
\, \dx\right)^{\frac{1}{p}},
\end{eqnarray}
by the definition of $\Psi$ at \eqref{defpsi} and Young's inequality.
At this point, using assumption  {{(F6)}} in the right hand side of \eqref{stimaquasifin}, we get 
\begin{eqnarray*}\label{stimaquasifin1}
	||V(Du)||_{L^{\infty}(B_{\bar\rho})} &\le& 
 C\left(\frac{ \|V(Du)\|_{\infty} }{\mathcal{L}\big(\|V(Du)\|_{\infty}\big)} + \frac{1}{( \bar R- \bar\rho )^{\frac{2n}{ p}}}\right)
\cdot \left(\int_{B_{ \bar R}}V(Du)^{p}
\, \dx\right)^{\frac{1}{p}}\cr\cr
&\le& C\left(\mathfrak{L}^{-1}\big(\|V(Du)\|_{\infty}\big) + \frac{1}{( \bar R- \bar\rho )^{\frac{2n}{ p}}}\right)
 \left(\int_{B_{ \bar R}}V(Du)^{p}
\, \dx\right)^{\frac{1}{p}},
    \end{eqnarray*}
    where we used Lemma \ref{GGPnew} and, to simplify the notation, we set $\mathfrak{L}(t)=\Big(\mathcal{L}^{-1}\Big)^*(t)$. Since  the function $\mathfrak{L}^{-1}(t)$ is  concave as it is the inverse of a polar function that is always convex, we may use 
    Young's inequality at \eqref{Young-psi}
 $$st\le \frac12\mathfrak{L}(s)+ \mathfrak{L}^*(2  t)$$
 %\colorbox{yellow}{perché $c$? In realtà $c=2$ perchè scriviamo $st=\frac12s\cdot2t$}
 This yields
\begin{eqnarray*}
	||V(Du)||_{L^{\infty}(B_{\bar\rho})} &\le&  \frac12||V(Du)||_{L^{\infty}(B_{ \bar R})}+ \mathfrak{L}^*\Bigg(C\left(\int_{B_{ \bar R}}V(Du)^{p}\, \dx\right)^{\frac{1}{p}}\Bigg)\cr\cr
&&+\frac{C}{( \bar R- \bar\rho )^{\frac{2n}{ p}}}\left(\int_{B_{ \bar R}}V(Du)^{p}
\, \dx\right)^{\frac{1}{p}}.
\end{eqnarray*}
Since previous estimate holds true for every $\rho_0<\bar\rho<\bar R<R_0$, the iteration Lemma \ref{lem:Giaq} implies the conclusion.
 \iffalse
 
 At this point, by \eqref{assumptionpsi} and \eqref{conseguenzaLemma2.1Psi}, we are allowed to use Young's inequality at \eqref{Young-psi}
 $$st\le \frac12\tilde\Psi^{-1}(s)+ \Big(\tilde\Psi^{-1}\Big)^*(c  t)$$
 where $\Big(\tilde\Psi^{-1}\Big)^*(t)$ denotes the Fenchel tranform of $\tilde\Psi^{-1}(s)$.
\tcr{In a similar way, by (G1) and \eqref{conseguenzaLemma2.1g}  
 we are allowed to use Young's inequality at \eqref{Fenchel}
 $$st\le \frac12\tilde g^{-1}(s)+ \Big(\tilde g^{-1}\Big)^*(c  t)$$
 where $\Big(\tilde g^{-1}\Big)^*(t)$ denotes the Fenchel tranform of $\tilde g^{-1}(s): = ??????$ }
 This yields
\begin{eqnarray*}
	||V(Du)||_{L^{\infty}(B_{\bar\rho})} &\le&  \frac12||V(Du)||_{L^{\infty}(B_{\bar R})}\cr\cr
 &&+ \Big(\tilde\Psi^{-1}\Big)^*\Bigg(\frac{C}{(\bar R-\bar \rho )^{\frac{n}{p}}}\left(\int_{B_{\bar R}}V(Du)^{p}
\, dx\right)^{\frac{1}{p}}\Bigg)\cr\cr
&& \tcr{+ \Big(\tilde g^{-1}\Big)^*\Bigg(\frac{C}{(\bar R-\bar \rho )^{\frac{n}{p}}}\left(\int_{B_{\bar R}}V(Du)^{p}
\, dx\right)^{\frac{1}{p}}\Bigg),}
\end{eqnarray*}
where, we used the pointwise inequality $|Du|\le V(Du)$.
Hence, we may use Lemma  \ref{lem:Giaq} to get
\begin{eqnarray}\label{apriori}
	||V(Du)||_{L^{\infty}(B_{\rho})} &\le&   \Big(\tilde\Psi^{-1}\Big)^* \Bigg(\frac{C}{( R- \rho )^{\frac{n}{p}}}\left(\int_{B_{ R}}V(Du)^{p}
\, dx\right)^{\frac{1}{p}}\Bigg)\cr\cr
&& \tcr{+ \Big(\tilde g^{-1}\Big)^*\Bigg(\frac{C}{(\bar R-\bar \rho )^{\frac{n}{p}}}\left(\int_{B_{\bar R}}V(Du)^{p}
\, dx\right)^{\frac{1}{p}}\Bigg),}
\end{eqnarray}
i.e. the conclusion.
\fi
 \end{proof}
Note that if $\mathcal{L}^{-1}$ is {proper,} lower semicontinuous and convex $(\mathcal{L}^{-1})^{**}=\mathcal{L}^{-1}$ {by Fenchel-Moreau theorem (for more details see for instance \cite{Aubin98}, Chapter 3, Theorem 3.1).}

\section{The a priori estimate in case $\mathbf{n=2}$}
Here, we complete the proof of the a priori estimate dealing with the case $n=2$.
 \begin{proof}[Proof of Theorem \ref{thmapriori} in the case $\mathbf{n=2}$ ] Let us fix a ball $B_{R_0}\Subset \Omega$, radii $\rho_{0}<\bar\rho<\rho<R<\bar R<R_{0}$
 and a cut-off function $\eta\in C^1_0(B_{R})$ such that $\eta=1$ on $B_{\rho}$ and $|D\eta|\le \frac{C}{R-\rho}$. Our starting point is the Sobolev inequality at Theorem \ref{sobtru} with the function $w=\eta(1+|Du|^2)^{\frac{p+\gamma}{4}}$ and $n=2$, i.e.
 \begin{eqnarray*}
&&{\Big|\Big|\eta(1+|Du|^2)^{\frac{p+\gamma}{4}}\Big|\Big|^2_{{\rm EXP}_2(B_R)}}\le \mathcal{C}_S\int_{B_R}\left|D\left(\eta (1+|Du|^2)^{\frac{p+\gamma}{4}}\right)\right|^2\,\dx\cr\cr
&\le& {C}
(p+\gamma)^2\int_{B_R}\eta^2 (1+|Du|^2)^{\frac{p+\gamma-2}{2}}|D^2 u|^2\, \dx+\frac{{C}}{(R-\rho)^2}\int_{B_R}(1+|Du|^2)^{\frac{p+\gamma}{2}}\,\dx.
 \end{eqnarray*}
 Using {\eqref{Caccioppoli}} to estimate the first term in the right hand side of previous inequality {and recalling \eqref{V(Du)}} 
 \begin{eqnarray*}
&&{\Big|\Big|\eta {V(Du)^{\frac{p+\gamma}{2}}}\Big|\Big|^2_{{\rm EXP}_2(B_R)}}\cr\cr
&\le& {C} (p+\gamma)^4   \, \int_{B_{R}}\eta^2  h^2 {V(Du)^{p + \gamma}} g^2\big(V(Du)\big) \, \dx\cr\cr
&&\quad+  \frac{{C} (p+\gamma)^2}{(R-\rho)^2}  \, \int_{B_{R}}  {V(Du)^{p + \gamma}}   {g\big(V(Du)\big)} \, \dx\cr\cr
&&\quad+\frac{{C}}{(R-\rho)^2}\int_{B_R}{V(Du)^{p + \gamma}} \,\dx\cr\cr
&\le& {C} {(p+\gamma)^4}   \, \int_{B_{R}}\eta^2  h^2\,V(Du)^{p + \gamma} g^2\big(V(Du)\big) \, \dx\cr\cr
&&\quad+  \frac{{C}(p+\gamma)^2}{(R-\rho)^2}  \, \int_{B_{R}}   V(Du)^{p + \gamma}  {g\big(V(Du)\big)} \, \dx,
 \end{eqnarray*}
 where we used that $g(1+t)\ge 1$, for every $t\ge 0$. 
 The a priori assumption $Du\in L^{\infty}_{\mathrm{loc}}(\Omega)$ and the monotonicity of the function $g(t)$ imply
 \begin{eqnarray}\label{primastiman=2}
&&\Big|\Big|\eta V(Du)^{\frac{p+\gamma}{2}}\Big|\Big|^2_{{\rm EXP}_2(B_R)}%=\Big|\Big|\eta(1+|Du|^2)^{\frac{p+\gamma}{4}}\Big|\Big|^2_{{\rm EXP}_2(B_R)}}
\cr\cr
&\le& {C}(p+\gamma)^4 g^2\big(||V(Du)||_{\infty}\big)  \, \int_{B_{R}}\eta^2  h^2 V(Du)^{p + \gamma}  \, \dx\cr\cr
&&\quad+  \frac{{C}(p+\gamma)^2}{(R-\rho)^2}g\big(||V(Du)||_{\infty}\big)  \, \int_{B_{R}}   V(Du)^{p + \gamma}  \, \dx\cr\cr
&\le& {C}(p+\gamma)^4 g^2\big(||V(Du)||_{\infty}\big)  \, \int_{B_{R}\cap\{{h}\le \varepsilon\}}\eta^2  h^2 V(Du)^{p + \gamma}  \, \dx\cr\cr
&&+ {C}(p+\gamma)^4 g^2\big(||V(Du)||_{\infty}\big)  \, \int_{B_{R}\cap\{{h}> \varepsilon\}}\eta^2  h^2 V(Du)^{p + \gamma}  \, \dx\cr\cr
&&\quad+  \frac{{C}(p+\gamma)^2}{(R-\rho)^2}g\big(||V(Du)||_{\infty}\big)  \, \int_{B_{R}}   V(Du)^{p + \gamma}  \, \dx\cr\cr
&\le& {C}\,\varepsilon^2(p+\gamma)^4 g^2\big(||V(Du)||_{\infty}\big)   \int_{B_{R}}\eta^2   V(Du)^{p + \gamma}  \, \dx\cr\cr
&&+{C}(p+\gamma)^4 g^2\big(||V(Du)||_{\infty}\big) \, \int_{B_{R}}\eta^2  h^2\chi_{\{{h}> \varepsilon\}} V(Du)^{p + \gamma}  \, \dx\cr\cr
&&\quad+  \frac{{C}(p+\gamma)^2}{(R-\rho)^2}g\big(||V(Du)||_{\infty}\big)  \, \int_{B_{R}}   V(Du)^{p + \gamma}  \, \dx,
\end{eqnarray}
where $\varepsilon>0$ will be chosen later.
Using H\"older's inequality at Proposition \ref{HoldOrlicz} to estimate the second integral in the right hand side of the previous estimate, we obtain
\begin{eqnarray}\label{primastiman=2}
&&{\Big|\Big|\eta V(Du)^{\frac{p+\gamma}{2}}\Big|\Big|^2_{{\rm EXP}_2(B_R)}}\cr\cr
&\le&  {C}\,\varepsilon^2(p+\gamma)^4 g^2\big(||V(Du)||_{\infty}\big)   \int_{B_{R}}\eta^2   V(Du)^{p + \gamma}  \, \dx\cr\cr
&&+ {C}(p+\gamma)^4 g^2\big(||V(Du)||_{\infty}\big)\big|\big|h^2\chi_{\{{h}> \varepsilon\}} \big|\big|_{L\log L(B_R)} {\Big|\Big|\eta^2 V(Du)^{{p+\gamma}}\Big|\Big|_{{\rm EXP}(B_R)}}\cr\cr
&&\quad+  \frac{{C}(p+\gamma)^2}{(R-\rho)^2}g\big(||V(Du)||_{\infty}\big)  \, \int_{B_{R}}   V(Du)^{p + \gamma}  \, \dx\cr\cr
&\le& {C}\,\varepsilon^2(p+\gamma)^4 g^2\big(||V(Du)||_{\infty}\big)   \int_{B_{R}}\eta^2   V(Du)^{p + \gamma}  \, \dx\cr\cr
&&\stackrel{\eqref{galquadrato}}{+} {C}(p+\gamma)^4 g^2\big(||V(Du)||_{\infty}\big)\big|\big|h^2\chi_{\{{h}> \varepsilon\}} \big|\big|_{L\log L(B_{R_0})} {\Big|\Big|\eta V(Du)^{^{\frac{p+\gamma}{2}}}\Big|\Big|^2_{{\rm EXP}_{2}(B_R)}}\cr\cr
&&\quad+  \frac{{C}(p+\gamma)^2}{(R-\rho)^2}g\big(||V(Du)||_{\infty}\big)  \, \int_{B_{R}}   V(Du)^{p + \gamma}  \, \dx.
%&\le& \frac12 \textcolor{violet}{\Big|\Big|\eta(1+|Du|^2)^{\frac{p+\gamma}{4}}\Big|\Big|^2_{{\rm EXP}_2(B_R)}}\cr\cr 
%&&\quad+ c(p+\gamma)^8 ||g\big(V(Du)\big)||^4_{\infty}||h^2||^2_{L\log^{\frac12}L(B_R)}\cr\cr
%&&\quad  +\frac{c(p+\gamma)^2}{(R-\rho)^2}||g\big(V(Du)\big)||_{\infty}  \, \int_{B_{R}}   %V(Du)^{p + \gamma}  \, dx,
\end{eqnarray}
 Now, using \eqref{alphallogl} and  that 
 $
\log(e + h^2) \le \, \log ((e + h)^2) = 2 \log(e + h), 
 $ 
 we get, {for any exponent $\alpha$ greater than 1}
\begin{eqnarray*}
 &&\big|\big|h^2\chi_{\{{h}> \varepsilon\}}\big|\big|_{L\log L(B_{R_0})}\le {C}(\alpha,{|B_{R_0}|})\left(\int_{B_{R_0}}h^2\chi_{\{{h}> \varepsilon\}}\log (e+h^2\chi_{\{{h}> \varepsilon\}})\,\dx\right)^{\frac{1}{\alpha}}\cr\cr
&&\qquad + \int_{B_{R_0}}h^2\chi_{\{{h}> \varepsilon\}}\log (e+h^2\chi_{\{{h}> \varepsilon\}})\,\dx\cr\cr
&=& {C}(\alpha,|B_{R_0}|) \left(\int_{B_{R_0}\cap\{{h}> \varepsilon\}}h^2\log (e+h^2)\,\dx\right)^{\frac{1}{\alpha}}+ \int_{B_{R_0}\cap\{{h}> \varepsilon\}}h^2\log (e+h^2)\,\dx
 \cr\cr  &\le & {C}(\alpha, c_\Phi,|B_{R_0}|)\left(\int_{B_{R_0}\cap\{{h}> \varepsilon\}}h^2\Phi^\mu(h)\,\dx\right)^{\frac{1}{\alpha}}+\int_{B_{R_0}\cap\{{h}> \varepsilon\}}h^2\Phi^\mu(h)\,\dx,
 \end{eqnarray*}
 where we used that, by virtue of assumption {(F5)$_2$}, 
$$\log(e+t)\le c_\Phi \Phi^\mu(t)\qquad\qquad \mu<\frac{1}{4}$$
{(notice that once more this estimate holds for $\Phi(h) > 1$ which is not restrictive to assume, otherwise our estimate turns to be trivial).}
 Hence, for the exponent $\theta$ appearing in assumption (F4), we have
 \begin{eqnarray*}
 &&\big|\big|h^2\chi_{\{{h}> \varepsilon\}}\big|\big|_{L\log L(B_{R_0})}\cr\cr 
 &\le& {C}(\alpha, c_\Phi,|B_{R_0}|)\left(\int_{B_{R_0}\cap\{{h}> \varepsilon\}}h^2\Phi^{{\theta}}(h)\frac{1}{\Phi^{{\theta-\mu}}(h)}\,\dx\right)^{\frac{1}{\alpha}} +\int_{B_{R_0}\cap\{{h}> \varepsilon\}}h^2\Phi^{{\theta}}(h)\frac{1}{\Phi^{{\theta-\mu}}(h)}\,\dx\cr\cr 
 &\le& \frac{C(\alpha, c_\Phi,|B_{R_0}|)}{\Phi^{\frac{\theta-\mu}{\alpha}}(\varepsilon)}\left(\int_{B_{R_0}\cap\{{h}> \varepsilon\}}h^2\Phi^{{\theta}}(h)\,\dx\right)^{\frac{1}{\alpha}}+\frac{1}{\Phi^{\theta-\mu}(\varepsilon)}\int_{B_{R_0}\cap\{{h}> \varepsilon\}}h^2\Phi^{{\theta}}(h)\,\dx,
\end{eqnarray*}
where we used the monotonicity of the function $\Phi$. Since by {\eqref{conseguenzasuPhi}} it holds 
\[
\Phi(\varepsilon) \ge \max \left\{\frac{1}{c_\Phi}, \Phi(1)\right\} =: K 
\]
and since $\alpha>1$, we have that $\left(\frac{\Phi(\varepsilon)}{{K}}\right)^{\frac{\theta-\mu}{\alpha}}\le \left(\frac{\Phi(\varepsilon)}{{K}}\right)^{\theta-\mu} $ that yields $\Phi^{\frac{\theta-\mu}{\alpha}}(\varepsilon)\le {{K}}^{\frac{(\theta-\mu)(1-\alpha)}{\alpha}}\Phi^{\theta-\mu}(\varepsilon)$. Hence
\begin{eqnarray*}
 &&\big|\big|h^2\chi_{\{{h}> \varepsilon\}}\big|\big|_{L\log L(B_{R_0})}\cr\cr 
 &\le& \frac{C(\alpha,\theta,\mu, {\Phi,} |B_{{R_0}}|)}{\Phi^{\frac{\theta-\mu}{\alpha}}(\varepsilon)}\left[\left(\int_{B_{R_0}\cap\{{h}> \varepsilon\}}h^2\Phi^{{\theta}}(h)\,\dx\right)^{\frac{1}{\alpha}}+\int_{B_{R_0}\cap\{{h}> \varepsilon\}}h^2\Phi^{{\theta}}(h)\,\dx\right].
\end{eqnarray*}
Now, since  $\theta>\mu+2$ by \eqref{deftheta}, we may choose $\alpha=\frac{\theta-\mu}{2}>1$ thus getting
\begin{eqnarray*}
 &&\big|\big|h^2\chi_{\{{h}> \varepsilon\}}\big|\big|_{L\log L(B_{R_0})}\le \frac{C(\theta,\mu, |B_{R_0}|)}{\Phi^{2}(\varepsilon)}\left[\left(\int_{B_{R_0}}h^2\Phi^{{\theta}}(h)\,\dx\right)^{\frac{2}{\theta-\mu}}+\int_{B_{R_0}}h^2\Phi^{{\theta}}(h)\,\dx\right].
\end{eqnarray*}
Using the previous estimate in \eqref{primastiman=2}, we obtain
\begin{eqnarray}\label{primastiman=2bis}
&&{\Big|\Big|\eta V(Du)^{\frac{p+\gamma}{2}}\Big|\Big|^2_{{\rm EXP}_2(B_R)}}\cr\cr
&\le&  {C}\varepsilon^2(p+\gamma)^4 g^2\big(||V(Du)||_{\infty}\big)   \int_{B_{R}}\eta^2   V(Du)^{p + \gamma}  \, \dx\cr\cr
&&+{\hat{C}}(p+\gamma)^4 g^2\big(||V(Du)||_{\infty}\big)\frac{{1}}{{\Phi^2(\varepsilon)}}\left[\left(\int_{{B_{R_0}}}h^2\Phi^{{\theta}}(h)\,dx\right)^{\frac{2}{\theta-\mu}}+\int_{{B_{R_0}}}h^2\Phi^{{\theta}}(h)\,\dx\right]\cr\cr
&&\qquad\qquad\cdot{\Big|\Big|\eta V(Du)^{^{\frac{p+\gamma}{2}}}\Big|\Big|^2_{{\rm EXP}_2(B_R)}}\cr\cr 
&&\quad+  \frac{{C}(p+\gamma)^2}{(R-\rho)^2}g\big(||V(Du)||_{\infty}\big)  \, \int_{B_{R}}   V(Du)^{p + \gamma}  \, \dx\cr\cr
&\le&  {C}\varepsilon^2(p+\gamma)^4 g^2\big(||V(Du)||_{\infty}\big)   \int_{B_{R}}\eta^2   V(Du)^{p + \gamma}  \, \dx\cr\cr
&&+{{\hat{C}}}(p+\gamma)^4 g^2\big(||V(Du)||_{\infty}\big)\frac{ {H_2^{4/(\theta-\mu)}}+H_2^2}{{\Phi^2(\varepsilon)}}{\Big|\Big|\eta V(Du)^{^{\frac{p+\gamma}{2}}}\Big|\Big|^2_{{\rm EXP}_2(B_R)}}\cr\cr 
&&\quad+  \frac{{C}(p+\gamma)^2}{(R-\rho)^2}g\big(||V(Du)||_{\infty}\big)  \, \int_{B_{R}}   V(Du)^{p + \gamma}  \, \dx,
\end{eqnarray}
where {the constant $\hat{C}$ depends also on $\Phi, \theta,\mu, |B_{R_0}|$ and where} we used as before the notation
\[
{H_2 := \left (\int_{B_{R_0}} h^2 \Phi^{\theta}(h) \, dx \right )^{1/2}.}
\]
Choosing now $\varepsilon>0$ such that
$${\hat{C}}(p+\gamma)^4 g^2\big(||V(Du)||_{\infty}\big)\frac{H_2^{4/(\theta-\mu)}+H_2^2}{{\Phi^2(\varepsilon)}}=\frac{1}{2} $$
$$\Longleftrightarrow\,\, \varepsilon=\Phi^{-1}\Big(\sqrt{2 {{\hat{C}}}}(p+\gamma)^{2} g\big(||V(Du)||_{\infty}\big)\sqrt{H_2^{4/(\theta-\mu)}+H_2^2}\Big),$$
estimate \eqref{primastiman=2bis} becomes
\begin{eqnarray}\label{primastiman=2ter}
&&{\Big|\Big|\eta V(Du)^{\frac{p+\gamma}{2}}\Big|\Big|^2_{{\rm EXP}_2(B_R)}}\le \frac{1}{2}{\Big|\Big|\eta V(Du)^{\frac{p+\gamma}{2}}\Big|\Big|^2_{{\rm EXP}_2(B_R)}}\cr\cr
&& + {C}\Big[(p+\gamma)^2 g\big(||V(Du)||_{\infty}\big)\Phi^{-1}\Big({\sqrt{2{\hat{C}}} (p+\gamma)^{2} g\big(||V(Du)||_{\infty}\big)\sqrt{H_2^{4/(\theta-\mu)}+H_2^2}}\Big)\Big]^2 \cr\cr &&\qquad\qquad \cdot \int_{B_{R}}\eta^2   V(Du)^{p + \gamma}  \, \dx\cr\cr
&&\quad + \frac{{C}(p+\gamma)^2}{(R-\rho)^2}g\big(||V(Du)||_{\infty}\big) \, \int_{B_{R}}   V(Du)^{p + \gamma}  \, \dx.
\end{eqnarray}
Reabsorbing the first integral in the right hand side by the left hand side, we infer
\begin{eqnarray}
&&  \frac{1}{2}{\Big|\Big|\eta V(Du)^{\frac{p+\gamma}{2}}\Big|\Big|^2_{{\rm EXP}_2(B_R)}}\cr\cr 
&\le & {C}\Big[(p+\gamma)^2 g\big(||V(Du)||_{\infty}\big)\Phi^{-1}\Big({\sqrt{2 {\hat{C}}} (p+\gamma)^{2} g\big(||V(Du)||_{\infty}\big)\sqrt{H_2^{4/(\theta-\mu)}+H_2^2}}\Big)\Big]^2 \cr\cr&&\qquad\qquad\cdot  \int_{B_{R}}\eta^2   V(Du)^{p + \gamma}  \, \dx\cr\cr
&&\quad  +\frac{{C}(p+\gamma)^2}{(R-\rho)^2}g\big(||V(Du)||_{\infty}\big)  \, \int_{B_{R}}   V(Du)^{p + \gamma}  \, \dx\cr\cr
&\le & {C}(p+\gamma)^4 g^{2}\big(||V(Du)||_{\infty}\big)\Phi^{-1}\Big(\mathcal{C}(p+\gamma)^{{4}}\Big)\Phi^{-1}\Big(g^{{2}}\big(||V(Du)||_{\infty}\big) \Big)   \int_{B_{R}}\eta^2   V(Du)^{p + \gamma}  \, \dx\cr\cr
&&\quad  +\frac{{C}(p+\gamma)^2}{(R-\rho)^2}g\big(||V(Du)||_{\infty}\big)  \, \int_{B_{R}}   V(Du)^{p + \gamma}  \, \dx,
\end{eqnarray}
where we used Lemma \ref{elem1} and $\mathcal{C}=\mathcal{C}(\lambda,\Lambda, H_2, {\Phi, \theta, \mu, |B_{R_0}|})$ is a constant independent of $\gamma$.
By Proposition \ref{normalpexp}, 
 recalling that $\eta=1$ on $B_\rho$, we get
 \begin{eqnarray}\label{terzastiman=2}
&&\left(\int_{B_{\rho}}V(|Du|)^{(p+\gamma)\frac{2^*}{2}} \,\dx\right)^{\frac{2}{\cdot2^*}}\cr\cr
&\le & {C} (p+\gamma)^4 g^{2}\big(||V(Du)||_{\infty}\big)\Phi^{-1}\Big(\mathcal{C}(p+\gamma)^{{4}}\Big)\Phi^{-1}\Big(g^{{2}}\big(||V(Du)||_{\infty}\big) \Big)    \int_{B_{R}}\eta^2   V(Du)^{p + \gamma}  \, \dx\cr\cr
&&\quad  +\frac{{C} (p+\gamma)^2}{(R-\rho)^2}g\big(||V(Du)||_{\infty}\big)  \, \int_{B_{R}}   V(Du)^{p + \gamma}  \, \dx, 
%&\le&  \frac{c(p+\gamma)^4 g\big(||V(Du)||_{\infty}\big)}{(R-\rho)^2}\Big[1+g\big(||V(Du)||_{\infty}\big)\Phi^{-1}\Big(\mathcal{C}(p+\gamma)^{{4}}\Big)\Phi^{-1}\Big(g^{{2}}\big(||V(Du)||_{\infty}\big) \Big)  \Big] \cr\cr&&\qquad\qquad\cdot  \int_{B_{R}}   V(Du)^{p + \gamma},
 \end{eqnarray}
where, %as before we supposed that $R-\rho<1$ and 
now {the constants depend also on $2^*$, which} denotes any finite exponent greater than $2$, since $n=2$. 
As in the case {$n\ge 3$}, defining the decreasing sequence of radii by setting
$$\rho_i=\bar \rho+\frac{ \bar R- \bar\rho}{2^i}$$
and the following increasing sequence of exponents
$$p_0=p \qquad\qquad {p_{i}}={p_{i-1}}\frac{2^*}{2}=p_0\left(\frac{2^*}{2}\right)^{i},$$
 we may iterate estimate \eqref{terzastiman=2}  on the concentric balls $B_{\rho_i}$  with exponents $p_i$,
thus obtaining
\begin{eqnarray}\label{stimapreitepasson=2}
&&\left ( \int_{B_{\rho_{i+1}}}  V(Du)^{p_{i+1}}\, \dx\right )^{\frac{1}{p_{i+1}}} 	\cr\cr
&\le&\!\! \displaystyle{\prod_{j=0}^{i}}\left({C} p_j^4g^2\big(||V(Du)||_{\infty}\big)\Phi^{-1}\Big(\mathcal{C}p_j^{{4}}\Big)\Phi^{-1}\Big(g^{{2}}\big(||V(Du)||_{\infty}\big) \Big)+\frac{Cp_j^{2}\,g\big(||V(Du)||_{\infty}\big)}{(\rho_j - \rho_{j+1})^{2}}     \right)^{\frac{1}{p_j}}\cr\cr&&\qquad\qquad\cdot\left(\int_{B_{\bar R}}V(Du)^{p}  \, \dx\right)^{\frac{1}{p}}\cr\cr
&=& \displaystyle{\prod_{j=0}^{i}}\left({C}  p_j^4g^2\big(||V(Du)||_{\infty}\big)\Phi^{-1}\Big(\mathcal{C}p_j^{{4}}\Big)\Phi^{-1}\Big(g^{{2}}\big(||V(Du)||_{\infty}\big) \Big)+\frac{C4^{j+1}p_j^{2}\,g\big(||V(Du)||_{\infty}\big)}{(\bar R - \bar\rho)^{2}}     \right)^{\frac{1}{p_j}}\cr\cr 
&&\qquad\qquad\cdot\left(\int_{B_{\bar R}}V(Du)^{p}  \, \dx\right)^{\frac{1}{p}}\cr\cr
&\le & \displaystyle{\prod_{j=0}^{i}}\left(g^2\big(||V(Du)||_{\infty}\big)\Phi^{-1}\Big(g^{{2}}\big(||V(Du)||_{\infty}\big)\Big)+\frac{g\big(||V(Du)||_{\infty}\big)}{(\bar R - \bar\rho)^{2}}     \right)^{\frac{1}{p_j}}\cr\cr 
&&\qquad \cdot\prod_{j=0}^{i}\left({C}  p_j^4\Phi^{-1}\Big(\mathcal{C}p_j^{{4}}\Big)+C4^{j+1}p_j^{2}\right)^{\frac{1}{p_j}}\cdot\left(\int_{B_{\bar R}}V(Du)^{p}  \, \dx\right)^{\frac{1}{p}}\cr\cr
&=:&\tilde\Theta_i\left(\int_{B_{\bar R}}V(Du)^{p}  \, \dx\right)^{\frac{1}{p}} .
\end{eqnarray}
Note that
$$\tilde\Theta_i \le \exp\left(\sum_{j=0}^i \frac{1}{p_j}\log \left( g^2\big(||V(Du)||_{\infty}\big)\Phi^{-1}\Big(g^{{2}}\big(||V(Du)||_{\infty}\big)\Big)+\frac{g\big(||V(Du)||_{\infty}\big)}{(\bar R - \bar\rho)^{2}}  \right) \right)$$
$$ \cdot \exp\left(\sum_{j=0}^i \frac{1}{p_j}\log\left({C}  p_j^4\Phi^{-1}\Big(\mathcal{C}p_j^{{4}}\Big)+C4^{j+1}p_j^{2}\right)\right)$$
$$\le \exp\left(\sum_{j=0}^\infty \frac{1}{p_j}\log\left( g^2\big(||V(Du)||_{\infty}\big)\Phi^{-1}\Big(g^{{2}}\big(||V(Du)||_{\infty}\big)\Big)+\frac{g\big(||V(Du)||_{\infty}\big)}{(\bar R - \bar\rho)^{2}}  \right) \right)$$
$$\cdot \exp\left(\sum_{j=0}^\infty \frac{1}{p_j}\log\left(1+\Phi^{-1}\Big(\mathcal{C}p_j^{{4}}\Big)  \right) \right)\cdot \exp\left(\sum_{j=0}^\infty \frac{1}{p_j}\log\left(1+C4^{j+1}  \right) \right)$$
$$ \cdot \exp\left(\sum_{j=0}^\infty \frac{1}{p_j}\log({C}  p_j^4)\right)=\tilde\Theta_1\cdot \tilde\Theta_2\cdot \tilde\Theta_3\cdot\tilde\Theta_4,$$
where, as before, we used the elementary inequality $\log(ab)\le \log(1+a)+\log(1+b)$, for $a,b>0$.
Since by the definition of $p_j$ we have $\sum_{j=0}^{{\infty}}\frac{1}{p_j}= \frac{2^*}{ p(2^*-2)}$,
hence
$$ \tilde\Theta_1 \le \left( g^2\big(||V(Du)||_{\infty}\big)\Phi^{-1}\Big(g^{{2}}\big(||V(Du)||_{\infty}\big)\Big)+\frac{g\big(||V(Du)||_{\infty}\big)}{(\bar R - \bar\rho)^{2}} \right)^{\frac{2^*}{ p(2^*-2)}}$$
as well as
$$\tilde\Theta_3\cdot \tilde\Theta_4  \le c(p,\lambda,\Lambda, 2^*).$$
By {(F5)$_2$} we derive
$$ \tilde\Theta_2\!\le \!\exp\left(\sum_{j=0}^{+\infty} \frac{1}{p_j}\Big( \mathcal{C} p_j^{{4}}\Big)^\mu\right)\le  C(p,2^*,\lambda,\Lambda, {H_2}, c_\Phi, {\theta, \mu, |B_{R_0}|}),$$
since, by our assumption, ${\mu<\frac14}$ {(as before this estimates holds when ${\mathcal{C} p_j^4 > 1}$ which is not restrictive to assume).}
Combining  the estimates above, we conclude that
\begin{eqnarray*}
\tilde \Theta_i&\le& C\left( g^2\big(||V(Du)||_{\infty}\big)\Phi^{-1}\Big(g^{{2}}\big(||V(Du)||_{\infty}\big)\Big)+\frac{g\big(||V(Du)||_{\infty}\big)}{(\bar R - \bar\rho)^{2}}\right)^{\frac{2^*}{ p(2^*-2)}}
\end{eqnarray*}
for every $i\in \mathbb{N}$, with $C=C(p,2^*,\lambda,\Lambda, {H_2}, {\Phi, \theta, \mu, |B_{R_0}|})$.
Therefore, we can let $i\to \infty$ in \eqref{stimapreitepasson=2} thus getting
\begin{eqnarray}
    ||V(Du)||_\infty&\le& C\left( g^2\big(||V(Du)||_{\infty}\big)\Phi^{-1}\Big(g^{{2}}\big(||V(Du)||_{\infty}\big)\Big)+\frac{g\big(||V(Du)||_{\infty}\big)}{(\bar R - \bar\rho)^{2}}  \right)^{\frac{2^*}{ p(2^*-2)}}\cr\cr 
&&\qquad\qquad\cdot\left(\int_{B_{\bar R}}V(Du)^{p}  \, \dx\right)^{\frac{1}{p}}\cr\cr 
&\le  & C\left (\Psi^{\frac{2^*}{p(2^*-2)}}\big(\|V(Du)\|_{\infty}\big) + g^{\frac{2\cdot2^*}{p(2^*-2)}}(\|V(Du)\|_{\infty})+ \frac{C }{(\bar R- \bar\rho )^{\frac{4\cdot 2^*}{p(2^*-2)}}} \right )
  \cr\cr 
&&\qquad\qquad\cdot \left(\int_{B_{ \bar R}}V(Du)^{p}
\, \dx\right)^{\frac{1}{p}},
\end{eqnarray}
by the definition of $\Psi$ and Young's inequality. This is precisely estimate \eqref{stimaquasifin} of previous subsection, and so the conclusion now comes arguing as in the case $n \ge 3.$
\end{proof}

 \section{The approximation}

\label{quattro}

In this section we give only a sketch of  the approximation procedure that comes along quite standard arguments (in view of the assumption $F(x, \xi) = \tilde{F}(x, |\xi|)$).
Let us start recalling an approximation Lemma whose proof can be achieved arguing exactly as in   \cite[Lemma 4.1]{CupGuiMas}.
\begin{lemma}\label{approx} Let $F:\Omega\times \mathbb{R}^n\to [0,+\infty)$ be a Caratheodory function satisfying assumptions {\rm (F1)--(F3)}. Then there exists a sequence of Caratheodory functions $F_k:\,\Omega\times \mathbb{R}^n\to [0,+\infty)$ monotonically convergent to $F(x,\xi)$ such that
\medskip

i) $F_k(x, \xi)=\tilde F_k(x, |\xi|)$, for a.e. $x\in \Omega$ and every $\xi\in \mathbb{R}^n$

ii) $F_k(x, \xi)\le  F_{k+1}(x, \xi)\le F(x,\xi)$, for a.e. $x\in \Omega$ , every $\xi\in \mathbb{R}^n$ and  $k\in\mathbb{N}$

iii) there exist $\ell$ independent of $k$ and $\ell_1$ depending also on $k$ such that
\[ \begin{cases}
  (1 + |\xi|^2)^{p/2}\le  F_k(x, \xi) \le \ell (1 + |\xi|^2)^{p/2} g({(1 + |\xi|^2)^{1/2}})
\cr\cr F_k(x, \xi) \le \ell_1(k) (1 + |\xi|^2)^{p/2}  
\end{cases}\eqno{\rm (A1)}\]
for a.e. $x\in \Omega$ , every $\xi\in \mathbb{R}^n$

iv) there exists a constant $\tilde\lambda$ depending only on $\lambda, p$ such that
\[\tilde\lambda \, (1 + |\xi|^2)^{\frac{p-2}{2}} |\mu|^2 \le {\sum_{i,j} (F_k)_{\xi_i \xi_j}(x, \xi) \mu_i \mu_j,}\eqno{\rm (A2)}\]
 for a.e. $x\in \Omega$ , every $\xi\in \mathbb{R}^n$

v)  there exist $\tilde\Lambda$ independent of $k$ and $\tilde\Lambda_1$ depending also on $k$ such that
\[
\left \{
\begin{array}{lll}
& |(F_k)_{{\xi_i \xi_j}}(x, \xi)| \le \, \tilde\Lambda \, (1 + |\xi|^2)^{\frac{p-2}{2}}  g({(1 + |\xi|^2)^{1/2}})\\[2mm] 
& |(F_k)_{{\xi_i \xi_j}}(x, \xi)| \le \, \tilde\Lambda_1(k) \, (1 + |\xi|^2)^{\frac{p-2}{2}}  
\end{array}
\right.
\eqno{\rm (A3)} 
\]
for a.e. $x\in \Omega$ , every $\xi\in \mathbb{R}^n$

vi)  there exist $C$ independent of $k$ and $C_1$ depending also on $k$ such that
\[
\left \{
\begin{array}{lll}
& |(F_k)_{\xi x}(x, \xi)| \le C h(x) (1 + |\xi|^2)^{\frac{p-1}{2}} g({(1 + |\xi|^2)^{1/2}})\\[2mm] 
&|(F_k)_{\xi x}(x, \xi)| \le C_1(k) h(x) (1 + |\xi|^2)^{\frac{p-1}{2}}  
\end{array}
\right.
\eqno{\rm (A4)}
\]
for a.e. $x\in \Omega$ , every $\xi\in \mathbb{R}^n$.
\end{lemma}
{In our case the original function $F$ possesses extra regularity (i.e. $\xi \mapsto F(x, \xi)$ is $\mathcal{C}^2(\real^n)$  and $x \mapsto {F_{\xi }(x, \xi)}$ is weakly differentiable, with $F_{\xi x}(x, \xi)$ a Carath\'eodory function) which transfers to the sequence of approximating functions $F_k.$}

Now for $u\in W^{1,p}(\Omega)$ a minimizer of the functional $\mathfrak{F}$, let $v_{k}\in u+W^{1,p}_0(B_{R_0})$ be the unique solution  to the problem
$$\min\left\{\mathfrak{F}_{k} (v,B_{R_0}):\,\,\, v\in u+W^{1,p}_0(B_{R_0})\right\},$$
where $$\mathfrak{F}_{k} (v,B_{R_0})=\int_{B_{R_0}}F_{k}(x, Dv)\,\dx.$$
By virtue of (A1)--(A4) and since $h(x)$ satisfies (F4) with a function $\Phi$ fulfilling (F5), we  may use \cite[Theorems 1 and 2]{CGHPdN20} %classical regularity results ({see for example \cite{Giusti}}) 
to deduce that  $v_{k}\in W^{1,\infty}_{\rm loc}(B_{R_0})\cap W^{2,{\min\{2,p\}}}_{\rm loc}(B_{R_0}) $  for every %$\varepsilon>0$ and 
$k\in \mathbb{N} $.
The left inequality in (A1) implies that
\begin{eqnarray*}
C_1 \int_{B_{R_0}}(1+|Dv_{k}|^2)^{\frac{p}{2}}\, \dx&\le& \int_{B_{R_0}}F_{k} (x,Dv_{k})\,\dx  \\[2mm]
&\le& \int_{B_{R_0}}F_{k} (x,Du)\,\dx,  \nonumber
\end{eqnarray*}
where we used the minimality of $v_{k}$.
%Since $F_{k} (x,Du)\in L^1(B_{R_0})$, for every $k\in \mathbb{N}$, taking the limit as $\varepsilon\to 0$, by the properties of the convolutions, we get
Hence, by  ii) of Lemma \ref{approx}, from previous estimate we deduce
\begin{equation}\label{normaLp}
\int_{B_{R_0}}(1+|Dv_{k}|^2)^{\frac{p}{2}}\, \dx \le    \int_{B_{R_0}}F_{k} (x,Du)\,\dx \le \int_{B_{R_0}}F (x,Du)\,\dx.    
\end{equation}
Therefore $\big(v_{k}\big)_k$ is a bounded sequence in $u+W^{1,p}(B_{R_0})$ 
\iffalse
Therefore $\big(v_{\varepsilon,k}\big)_\varepsilon$ is a bounded sequence in $u+W^{1,p}(B_{R_0})$ and so, for every $k\in \mathbb{N}$, there exists $v_k\in u+W^{1,p}(B_{R_0})$ such that \begin{equation}\label{weakLp}
  v_{\varepsilon,k} \rightharpoonup v_k \qquad \text{in}\,\, u+W^{1,p}(B_{R_0}).  
\end{equation}
The lower semicontinuity of the $L^p$ norm implies
\[
\int_{B_{R_0}}(1+|Dv_{k}|^2)^{\frac{p}{2}}\, \dx\le \liminf_{\varepsilon\to 0}\int_{B_{R_0}}(1+|Dv_{\varepsilon,k}|^2)^{\frac{p}{2}}\, \dx \le \int_{B_{R_0}}F (x,Du)\,\dx 
\]
\fi
and so 
there exists $v\in u+ W^{1,p}(B_{R_0})$ such that \begin{equation}\label{weakLpbis}
  v_{k} \rightharpoonup v \qquad \text{in}\,\, u+W^{1,p}(B_{R_0}).  
\end{equation}
On the other hand, since $v_{k}\in W^{1,\infty}_{\rm loc}(B_{R_0}) \cap W^{2,{\min\{2,p\}}}_{\rm loc}(B_{R_0})$, we can use the a priori estimate at \eqref{apriori} to deduce that
\begin{eqnarray}\label{apriori2}
&&||1+Dv_{k}||_{L^\infty \left (B_{\frac{R_0}{2}}\right )}\le  \Big( \mathcal{L}^{-1}\Big)^{**}\left (C\int_{B_{R_0}}  \, (1 + |Dv_{k}|^2)^{\frac{p}{2}} \dx \right) ^{\frac{2}{p}}\cr\cr
&&\qquad\qquad\qquad +\frac{C}{ R_0^{\frac{4\cdot2^*}{p(2^*-2)}}}\left(\int_{B_{  R_0}}(1 + |Dv_{k}|^2)^{\frac{p}{2}}
\, \dx\right)^{\frac{1}{p}} \cr\cr
&\le&  \Big( \mathcal{L}^{-1}\Big)^{**}\left(C\int_{B_{R_0}}F (x,Du)\,dx\right)\cr\cr
&&\qquad\qquad\qquad +\frac{C}{ R_0^{\frac{4\cdot2^*}{p(2^*-2)}}}\left(\int_{B_{  R_0}}F (x,Du)
\, \dx\right)^{\frac{1}{p}},
\end{eqnarray}
 where in the last inequality we used \eqref{normaLp} and the monotonicity of the function $\Big( \mathcal{L}^{-1}\Big)^{**}(t)$. 
\\
Previous inequality yields that
\begin{equation}\label{weakLpbi}
  v_{\varepsilon,k} \stackrel{*}\rightharpoonup v \qquad \text{weakly * in}\,\, W^{1,\infty}_{\mathrm{loc}}(B_{R_0})  
\end{equation}
and, by the lower semicontinuity of the norm, also
\begin{eqnarray}\label{apriori2bis}
||1+Dv||_{L^\infty \left(B_{\frac{R_0}{2}}\right )}
&\le& \Big( \mathcal{L}^{-1}\Big)^{**}\left(\int_{B_{R_0}}F (x,Du)\,\dx\right).
\end{eqnarray}
At this point, the conclusion of the proof, that is showing that $v=u$ a.e. in $B_{R_0}$, follows by standard arguments that can be found for instance in \cite{EMM16,CGGP,CGHPdN20}.

\vspace{10mm}

\noindent {\bf Acknowledgments.}
M. Eleuteri and A. Passarelli di Napoli are members of the Gruppo
Nazionale per l’Analisi Matematica, la Probabilità e le loro Applicazioni (GNAMPA) of the Istituto
Nazionale di Alta Matematica (INdAM). \smallskip
\par\noindent
{\bf Funding}. This research was partly funded by:
\\ (i)   GNAMPA Project 2026, grant number E53C25002010001 ``Esistenza e regolarita' per soluzioni di equazioni ellittiche e paraboliche anisotrope''(M. Eleuteri, A. Passarelli di Napoli).
\smallskip
\par\noindent
{\bf AI declaration:} no AI has been used to produce this research.

\bigskip
\noindent\footnotesize{\textsc{Michela Eleuteri , Dipartimento di Scienze Fisiche, Informatiche e Matematiche, 
Universit\`{a} di Modena e Reggio Emilia, via Campi 213/b - 41121 Modena, Italy}}

{\it Email address:}\,\,{\texttt{michela.eleuteri@unimore.it}}
\bigskip

\noindent\footnotesize{\textsc{Antonia Passarelli di Napoli, Dipartimento di Matematica e Applicazioni "R.
Caccioppoli", Universit\`{a} di Napoli ``Federico
II", via Cintia - 80126 Napoli, Italy}}

{\it Email address:}\,\,{\texttt{antpassa@unina.it}}
\end{document}